\documentclass{article}
\usepackage{amssymb}
\usepackage{dsfont}
\usepackage{textcomp}
\usepackage{amsmath}
\usepackage{pgf}
\usepackage[colorlinks=true]{hyperref}
\usepackage{txfonts}
\newtheorem{theo}{Theorem}[section]
\newtheorem{definition}[theo]{Definition}
\newenvironment{defi}{\begin{definition}\rm}{\end{definition}}
\newtheorem{remarque}[theo]{Remark}
\newenvironment{remark}{\begin{remarque}\rm}{\end{remarque}}
\newtheorem{exemple}[theo]{Example}
\newenvironment{ex}{\begin{exemple}\rm}{\end{exemple}}
\newtheorem{lemma}[theo]{Lemma}
\newtheorem{propo}[theo]{Proposition}
\newtheorem{coro}[theo]{Corollary}
\newtheorem{claim}[theo]{Claim}
\newtheorem{nota}[theo]{Notation}
\newenvironment{notation}{\begin{nota}\rm}{\end{nota}}

\newtheorem{prf}{\it{Proof.}}

\newenvironment{demonstration}{\begin{prf}\rm}{\hfill$\Box$\end{prf}}

\newcommand{\n}{\par\noindent}

\newcommand{\sn}{\par\smallskip\noindent}
\newcommand{\mn}{\par\medskip\noindent}
\newcommand{\bn}{\par\bigskip\noindent}

\title{Hardy type derivations on fields of exponential logarithmic series}
\author{Salma Kuhlmann.\\ Universit\"at Konstanz,\\ Fachbereich Mathematik und Statistik\\
78457 Konstanz, Germany.\\
Email: salma.kuhlmann@uni-konstanz.de\\
\bn Webpage: http://www.math.uni-konstanz.de/~kuhlmann/\\
Micka\"{e}l Matusinski 
\footnote{Current address : Institut de Math\'ematiques de Bordeaux, Universit\'e Bordeaux 1, 351 cours de la Lib\'eration - F 33405 Talence cedex, France}.\\ Universit\"at Konstanz,\\ Fachbereich Mathematik und Statistik\\ 78457 Konstanz, Germany.\\
Email: mickael.matusinski@u-bordeaux1.fr\\
Webpage: http://sites.google.com/site/mickaelmatusinski/}
\begin{document}
\maketitle
\begin{abstract}    % type your abstract below
We consider the valued field $\mathds{K}:=\mathbb{R}((\Gamma))$ of
formal series (with real coefficients and monomials in a
totally ordered multiplicative group $\Gamma\>$). We investigate how
to endow $\mathds{K}$ with a logarithm $l$, which
satisfies some natural properties such as commuting with infinite
products of monomials. In
\cite{matu-kuhlm:hardy-deriv-gener-series}, we studied derivations
on $\mathds{K}$. Here, we investigate compatibility conditions
between the logarithm and the derivation, i.e. when the logarithmic
derivative is the derivative of the logarithm. We analyse sufficient conditions on a given derivation to construct a compatible logarithm via integration of logarithmic derivatives. In \cite{kuhl:ord-exp}, the first author described the exponential closure $\mathds{K}^{\rm{ EL }}$ of $(\mathds{K},l)$. Here we show
how to extend such a log-compatible derivation on $\mathds{K}$ to
$\mathds{K}^{\rm{ EL }}$.
\end{abstract}

Keywords: generalized series field, pre-logarithm, derivation,
valuation, exponential closure.

Classification msc2000: primary 12J10,
12J15, 12L12, 13A18; secondary: 03C60, 12F05,
12F10, 12F20

\tableofcontents{}

\section{Introduction}
 Consider the valued
field $\mathds{K}:=\mathbb{R}((\Gamma))$ of generalized series, with
real coefficients and monomials in a totally ordered multiplicative
group $\Gamma\>$. We undertook  the investigation of these fields in a series of publications
\cite{kuhl:exp-pow-series},  \cite{kuhl:ord-exp},
\cite{kuhlm-fornasiero} and
\cite{matu-kuhlm:hardy-deriv-gener-series}. We endeavor to endow these formal algebraic objects with the analogous of classical analytic structures, such as exponential and logarithmic maps, derivation, integration and difference operators. Hardy fields, extensively studied by M. Rosenlicht, are the natural domain for asymptotic analysis. Our investigations thus lead us to analyse the relationship between Hardy fields and generalized series fields.  This paper is a further step in this direction. In particular, we interprete here some key ideas of \cite{rosenlicht:rank} in the formal setting of generalized series.
 \mn In \cite{kuhl:exp-pow-series}, we proved that
if $\Gamma \ne 0$, then $\mathds{K}$ cannot be endowed with a
logarithm (i.e. an isomorphism of ordered groups from its
multiplicative group of positive elements onto its additive group).
We established however that $\mathds{K}$ always admits a
pre-logarithm, i.e. a non surjective logarithm. In this paper, we
take a closer look at this aspect. We investigate how to endow
$\mathds{K}$ with a (non surjective) logarithm $l$, which moreover
satisfies some natural properties such as commuting with infinite
products of monomials. \mn
 In \cite{matu-kuhlm:hardy-deriv-gener-series}, we studied derivations on $\mathds{K}$ and introduced in particular Hardy type derivations (that is, derivations that behaves like derivations in a Hardy field). For the analysis of the derivations on $\mathds{K}$, we worked with the chain of fundamental monomials $(\Phi,\preccurlyeq)$ of $\Gamma$ (see Section \ref{sect:defi}). We gave a necessary and sufficient condition for a map $d\ :\ \Phi\rightarrow\mathds{K}$ to extend naturally to such a derivation. Here, we investigate compatibility conditions between the logarithm and the derivation, i.e. when the logarithmic derivative is the derivative of the logarithm. \mn In
\cite{kuhl:ord-exp}, the first author described the exponential
closure $\mathds{K}^{\rm{ EL }}$ of $(\mathds{K},l)$. Here we show
how to extend such a log-compatible derivation on $\mathds{K}$ to
$\mathds{K}^{\rm{ EL }}$. This exponential closure $\mathds{K}^{\rm{
EL }}$ is an infinite towering extension, starting with a
pre-logarithmic series field, i.e. a generalized series field
endowed with a pre-logarithm (see Definition \ref{defi:prelog}). Thus we begin in Section \ref{sect:prelog} by proving a criterion for a derivation on a pre-logarithmic series field to be compatible
 (see Proposition \ref{propo:compat}). This result is applied in Section \ref{section:hardy-derivation}. There, the main Theorem \ref{theo:existence-prelog} deals with a Hardy type series derivation $d$, and gives sufficient conditions on $d$ to define a $d$-compatible pre-logarithm. This pre-logarithm is constructed by a process of ``iterated asymptotic integration'' of the logarithmic derivatives (Lemma \ref{next}). This process is based on the computation of specific asymptotic integrals, which we do in Section \ref{section:asymp-integ}. This allows us to provide many examples in Section \ref{section:decr-autom}. In Section \ref{sect:EL}, given some pre-logarithmic series field endowed with a Hardy type derivation, we show how to extend it to the corresponding exponential closure. Note that this has been considered for fields of transseries in \cite[Ch. 4.1.4]{Schm01}. However, our pre-logarithmic field $(\mathds{K},l)$ does not necessarily satisfy Axiom (T4)
  of \cite[Definition 2.2.1]{Schm01}. The last Section \ref{sect:integr} is devoted to the questions of asymptotic integration and integration on EL-series fields.
    \mn
In forthcoming papers, we extend our investigations to study Hardy
type derivations on the field of Surreal Numbers \cite{conway_numb-games}, and investigate difference operators on generalized series fields.

%-----------------------------------------------------------------------------------------------------------------
\section{Preliminaries.}\label{sect:defi}

We summarize notation and terminology from \cite{matu-kuhlm:hardy-deriv-gener-series}.
Recall the following corollary to Ramsey's theorem \cite{rosen:lin-ord}:
\begin{lemma}\label{lemme:suite}
Let $\Gamma$ be a totally ordered set. Every sequence
$(\gamma_n)_{n\in\mathbb{N}}\subset \Gamma$ has an infinite
sub-sequence which is either constant, or strictly increasing, or
strictly decreasing.
\end{lemma}

\subsection{Hahn groups.}
\begin{defi}\label{defn: hahn group} Let $(\Phi,\preccurlyeq)$ be a totally ordered set,
the set of \textbf{fundamental monomials}. Consider the set $\mathrm{\textbf{H}}(\Phi)$
of formal products $\gamma$ of the form
\[\gamma=\prod_{\phi\in\mbox{ supp }\gamma}\phi^{\gamma_\phi}\]
 where $\gamma_\phi \in\mathbb{R}$, and \textbf{support} of
 $\gamma$, $\textrm{supp}\ \gamma := \{\phi\in\Phi\ |\ \gamma_{\phi}\neq 0\}$,
is an anti-well-ordered subset of $\Phi$. We will refer to $\gamma_\phi$ as the \textbf{exponent} of $\phi$.
Multiplication of formal
products is pointwise, and $\mathrm{\textbf{H}}(\Phi)$ is an abelian group with identity $1$.
We endow $\mathrm{\textbf{H}}(\Phi)$ with the anti lexicographic ordering $\preccurlyeq$ which
extends $\preccurlyeq$ of
 $\Phi$. Note that
$\phi \succ 1$ for all $\phi\in \Phi$. The totally
ordered abelian group  $\mathrm{\textbf{H}}(\Phi)$ is the \textbf{Hahn group} over $\Phi$, which elements are called the \textbf{(generalized) monomials}.
%We denote by $\mathrm{\textbf{H}}_{\mathrm{fin}}(\Phi)$ its subgroup of monomials with finite support.
The set $\Phi$ is the \textbf{rank}. By Hahn's embedding theorem
\cite{hahn:nichtarchim}, every ordered abelian group $\Gamma$ with rank $\Phi$ %(and archimedean components $\mathbb{R}$; \cite{kuhl:ord-exp})
can be seen as a subgroup of $\mathrm{\textbf{H}}(\Phi)$. 
\sn \textit{From now on, we fix a totally ordered set $(\Phi,\preccurlyeq)$ and a subgroup $\Gamma$ of $\mathrm{\textbf{H}}(\Phi)$}.
 \end{defi}

 \begin{defi}\label{defn:leading fundamental}
The \textbf{leading fundamental monomial} of $1\ne \gamma\in \Gamma$ is $\textrm{LF}(\gamma):=\max(\textrm{supp}\
 \gamma)\>$, and $\mbox{LF}(1):=1$.  This map
 verifies the \textbf{ultrametric triangular inequality} :
\[\forall \alpha,\beta\in\Gamma,
\mbox{ LF }(\alpha\beta)\preccurlyeq\max\{\mbox{ LF }(\alpha),\mbox{
LF }(\beta)\}\]. \sn The \textbf{leading exponent} of $1\ne \gamma
\in \Gamma$ is the exponent of $\mbox{ LF }(\gamma)$. We denote it
by $\mbox{ LE }(\gamma)$. For $\alpha \in \Gamma$ set $|\alpha|: =
\mbox{ max } \{\alpha\>, \> \alpha ^{-1}\}$.
\end{defi}

%------------------------------------------------------
\subsection{generalized series fields.}
Below, we adopt our notation as in \cite{matu-kuhlm:hardy-deriv-gener-series}.
\begin{defi}\label{defi:generalized series}
Throughout this paper,
$\mathds{K}=\mathbb{R}\left(\left(\Gamma\right)\right)$ will denote
the \textbf{generalized series field}. As usual, we write these series
$a=\displaystyle\sum_{\alpha\in\textrm{Supp}\ a}a_\alpha \alpha$,
and denote by $0$ the series with
 empty support. Here $\textrm{Supp}\ a=\{\alpha\in\Gamma\ |\ a_\alpha\neq 0\}$ is anti-well-ordered in $\Gamma$.
\sn For $a\in\mathds{K}^*$, its \textbf{leading monomial} is:
$\mbox{ LM }(a) := \max\left(\textrm{Supp}\ a\right)\in\Gamma\>$.
The map $\mbox{ LM }\>:\> \mathds{K}^* \rightarrow \Gamma$ is the  \textbf{canonical valuation} on $\mathds{K}$.
The \textbf{leading coefficient} of $a$ is $\mbox{ LC
}(a):= a_{\textrm{ LM }(a)}\in \mathbb{R}$. For nonzero
$a\in\mathds{K}$, the term $\mbox{ LC }(a)\mbox{ LM }(a)$ is called
the \textbf{leading term} of $a$, that we denote $\mbox{LT}(a)$.
We extend the notions of \textbf{leading fundamental
monomial} and of \textbf{leading exponent} to $\mathds{K}^*$ by setting $\mbox{ LF }(a):=\mbox{ LF }(\mbox{ LM }(a))$, respectively
 $\mbox{ LE }(a):=\mbox{ LE }(\mbox{ LM }(a))$.
\mn  We extend the ordering $\preccurlyeq$ on $\Gamma$ to a {\bf
dominance relation} on $\mathds{K}$ by setting $a\preccurlyeq
b\Leftrightarrow \mbox{ LM }(a)\preccurlyeq \mbox{ LM }(b)$. We
write: $a\asymp b \Leftrightarrow  \mbox{ LM }(a)= \mbox{ LM }(b)$,
and: $a\sim b \Leftrightarrow  \mbox{ LT }(a)= \mbox{ LT }(b)$. Let
$a\succ 1,\ b\succ 1$ be two elements of $\mathds{K}$. $a$ and $b$
are \textbf{comparable} if and only if $\mbox{ LF }(a)=\mbox{ LF
}(b)$. We also set $|a|: = |\mbox{ LM }(a)|$.
  \mn The \textbf{anti lexicographic ordering} on $\mathds{K}$ is defined as follows:
  $\forall a\in\mathds{K},\ a\leq
0\Leftrightarrow \mbox{ LC }(a)\leq 0\>$. We denote as usual
$\mathds{K}^*:=\mathds{K}\setminus\{0\}$, and
$\mathds{K}_{>0}:=\{a\in\mathds{K}\ |\ a>0\}$. Note that
$(\mathds{K}_{>0},\cdot)$ is an ordered abelian group. \end{defi}

\begin{remark}
The results in this paper hold for the generalized series field with
coefficients in an arbitrary ordered exponential field $\mathcal{C}$ \cite{kuhl:ord-exp} containing
$\mathbb{R}$ (instead of $\mathbb{R}$).
\end{remark}

%--------------------------------------------------
\subsection{Pre-logarithmic sections.}

\begin{defi}
We denote by $\mathds{K}^{\preccurlyeq 1} := \{a\in\mathds{K}\ |\
a\preccurlyeq 1\}$ the \textbf{valuation ring} of $\mathds{K}$.
Similarly, we denote by $\mathds{K}^{\prec 1} := \{a\in\mathds{K}\
|\ a\prec 1\}$ the \textbf{maximal ideal} of
$\mathds{K}^{\preccurlyeq 1}$. We have $\mathds{K}^{\preccurlyeq 1}=
\mathbb{R}\oplus \mathds{K}^{\prec 1}$. We denote by
$\mathds{K}^{\succ 1} := \mathbb{R}\left(\left(\Gamma^{\succ
1}\right)\right)$, the \textbf{subring of purely infinite series}.
\end{defi}
We will use repeatedly the following direct sum, respectively direct product, decompositions
of the ordered abelian groups $(\mathds{K},+,\leq)$, respectively $(\mathds{K}_{>0},\cdot,\leq)$
 \cite[Ch. 1]{kuhl:ord-exp}: \begin{center}
$\begin{array}{lclcccc}
\mathds{K} &=& \mathds{K}^{\succ 1}&\oplus &\mathbb{R} &\oplus& \mathds{K}^{\prec 1}\\
\mathds{K}_{>0} &=& \Gamma & .&\ \mathbb{R}_{>0} &.&\  (1 +\mathds{K}^{\prec 1})
\end{array}$
\end{center}

%The first decomposition is the classical one associated to the natural valuation \textrm{LM}. The second one means that any positive element $a$ can be written $a=A_\alpha\alpha (1+b)$ where $\alpha=\textrm{LM}(a)$, $A_\alpha=LC(a)$ and $h_a\in\mathds{K}^{\prec 1}$ (indeed, it suffices to divide $a$ by its leading term $A_\alpha\alpha$).

\begin{defi}\label{defi:prelog} Let $\mathds{K}$ be a field of generalized series.
\begin{itemize}
\item The natural \textbf{logarithm on 1-units} is the following isomorphism of ordered groups
 \cite{alling:exp_closed_fields}:
\begin{equation}\label{eq:log-1-unit}
\begin{array}{llcl}
l_1\ :\ &1+\mathds{K}^{\prec 1}& \rightarrow&\mathds{K}^{\prec 1}\\
&1+\epsilon&\mapsto&\displaystyle\sum_{n\geq 1}(-1)^{n-1}\displaystyle\frac{\epsilon^n}{n}
\end{array}
\end{equation}
    \item A \textbf{pre-logarithmic section} $l$ of $\mathds{K}$ is an embedding of ordered groups
\begin{center}
$l : (\Gamma , \cdot, \preccurlyeq) \rightarrow (\mathds{K}^{\succ 1}, +, \leq)$.
\end{center}
\item The \textbf{pre-logarithm} on $\mathds{K}$ induced by a pre-logarithmic section $l$ is the
embedding of ordered groups defined by:
\begin{center}
$\begin{array}{llcl}
l\ :\ &(\mathds{K}^{>0},\cdot,\leq)&\rightarrow&(\mathds{K},+,\leq)\\
&a=a_\alpha\alpha(1+\epsilon_a)&\mapsto &l(a):=\log(a_\alpha)+l(\alpha)+l_1(1+\epsilon_a)
\end{array}$
\end{center}
\sn where $\log$ is the usual logarithm on positive real numbers and
$l_1$ the logarithm on 1-units. The pair $(\mathds{K},l)$ is then
called a \textbf{pre-logarithmic series field}.
\end{itemize}
\end{defi}

In particular, we are interested in pre-logarithmic sections that verify the \textbf{Growth Axiom Scheme} :

\begin{defi}\label{defi:GA} Let $(\mathds{K},l)$ a pre-logarithmic series field.
We say that the Growth Axiom Scheme holds if and only if we have :
\begin{description}
\item[(GA)] $\forall \alpha\in\Gamma^{\succ 1}$, $l(\alpha)\prec \alpha$.
\end{description}
\end{defi}

Axiom (GA) is satisfied by non archimedean models of the theory $\mbox{Th}(\mathbb{R},\exp)$
of the ordered field of real numbers with the exponential function (see \cite[Ch. 3]{kuhl:ord-exp}
for more details).

%-------------------------------------------------------------------------------------

\section{Pre-logarithms and Derivations.}\label{sect:prelog}
\subsection{Defining pre-logarithms on generalized series fields.}

We consider \textsl{pre-logarithmic sections} on a generalized series field $\mathds{K}$
(see Definition \ref{defi:prelog}), which satisfy the following property:

\begin{defi}\label{defi:formule_pre-log}
A map $l:\Gamma\rightarrow \mathds{K}^{\succ 1}$ is a \textbf{series morphism} if it satisfies the following axiom:
\begin{description}
\item[(L)] $\forall\alpha=\displaystyle\prod_{\phi\in\textrm{supp}\ \alpha}\phi^{\alpha_\phi}\in\Gamma,
\ \ l(\alpha)=\displaystyle\sum_{\phi\in\textrm{supp}\ \alpha}\alpha_\phi l(\phi)$.
\end{description}
\end{defi}
Note that a series morphism $l$ is in particular a group homomorphism.
Moreover, a series morphism $l$ is a pre-logarithmic section if and only if it is order preserving
(i.e. for any $\alpha \prec \beta$ in $\Gamma$, we have $l(\alpha)< l(\beta)$).

We are interested in the following setting: given a map
 $l_\Phi:\Phi\rightarrow \mathds{K}^{\succ 1}\setminus\{0\}$, we study necessary and sufficient conditions
 so that $l_\Phi$ extends to a series morphism $l_\Gamma:\Gamma\rightarrow \mathds{K}^{\succ 1}$.

Recall the following definition from \cite{matu-kuhlm:hardy-deriv-gener-series}:
\begin{defi}\label{defi:summ_fam}
Let $I$ be an infinite index set and $\mathcal{F}=(a_i)_{i\in I}$ be a family
of series in $\mathds{K}$. Then $\mathcal{F}$ is said to be
\textbf{summable} if the two following properties hold:
\begin{description}
    \item[(SF1)]  $\textrm{Supp}\ \mathcal{F}:=\displaystyle\bigcup_{i\in I}\textrm{Supp}\ a_i$ (the support of the family) is
    an anti-well-ordered subset of $\Gamma$.
\item[(SF2)] For any $\alpha\in\textrm{Supp}\ \mathcal{F}$, the set $S_\alpha:=\{i\in I\ |\ \alpha\in\textrm{Supp}\ a_i\}\subseteq I$ is
finite.
\end{description}
Write $a_i=\displaystyle\sum_{\alpha\in\Gamma}a_{i,\alpha}\alpha$,
and assume that $\mathcal{F}=(a_i)_{i\in I}$ is summable. Then
\begin{center}
$\displaystyle\sum_{i\in I} a_i:=\displaystyle\sum_{\alpha\in
\textrm{Supp}\ \mathcal{F}}\displaystyle\>(\sum_{i\in
S_\alpha}a_{i,\alpha}\>)\alpha\in\mathds{K}$
\end{center}
is a well defined element of $\mathds{K}$ that we call the \textbf{sum} of $\mathcal{F}$.
\end{defi}

\begin{defi}\label{defi:sommabilite}
Let
\begin{center}$\begin{array}{lccl}
l_\Phi\ :&\Phi&\rightarrow& \mathds{K}^{\succ 1}\backslash\{0\}\\
&\phi &\mapsto& l_\Phi(\phi)
\end{array}$
\end{center}
be a map. We say that $l_\Phi$ \textbf{extends to a series morphism on $\Gamma$} if the following property holds:
\begin{description}
\item[(SL)] For any anti-well-ordered subset $E\subset\Phi$,
the family $\left(l_\Phi(\phi)\right)_{\phi\in
E}$ is summable. \end{description} Then the series
morphism $l_\Gamma$ on $\Gamma$ \textbf{induced} by $l_\Phi$ is defined
to be the map \[l_\Gamma :\ \Gamma\rightarrow\mathds{K}^{\succ 1}\] obtained
through the axiom (L) (which clearly makes sense by (SL)).
\end{defi}

Note that, if the series morphism $l_\Gamma$ is a pre-logarithmic section, then it extends to a pre-logarithm $l$ on $\mathds{K}_{>0}$ as in Definition \ref{defi:prelog}. We are interested moreover in pre-logarithms which verify (GA) (see Definition \ref{defi:GA}).

 In the next Proposition \ref{propo:pre-log}, we provide a necessary and sufficient condition on a
 map $l_\Phi :\Phi\rightarrow \mathds{K}$ so that the properties (SL) and (GA) hold. (In the sequel, we drop the subscripts $\Phi$ and $\Gamma$ of $l_\Phi$
 and $l_\Gamma$ to relax the notation).

\begin{propo}\label{propo:pre-log}
A map $l : \Phi\rightarrow\mathds{K}^{\succ 1}\setminus\{0\}$ extends to a series morphism on $\Gamma$ if and only if the following condition fails:
\begin{description}
    \item[(HL1)] there exists a strictly decreasing sequence $(\phi_n)_{n\in\mathbb{N}}\subset\Phi$ and an increasing sequence
 $(\lambda^{(n)})_{n\in\mathbb{N}}\subset\Gamma$ such that for any $n$, $\lambda^{(n)}\in\textrm{Supp}\ l(\phi_n)$.
\end{description}
Moreover, such an extension $l$ is a pre-logarithmic section if and only if we have:
\begin{description}
    \item[(HL2)] $l$ is an embedding of ordered sets, i.e. for any $\phi\prec\psi \in\Phi$, $0<l(\phi)< l(\psi)$.
%for any $\phi\prec\psi\in\Phi$, for any $\lambda^{(\phi)}\in\textrm{Supp}\ l(\phi)$,
%$\lambda^{(\psi)}\in\textrm{Supp}\ l(\psi)$, we have
% $LF\left(\displaystyle\frac{\lambda^{(\phi)}}{\lambda^{(\psi)}}\right)\prec\phi$.
\end{description}
Moreover, such a pre-logarithmic section $l$ satisfies (GA) if and only if we have:
\begin{description}
    \item[(HL3)] for any $\phi\in\Phi$, $\mbox{ LF }(l(\phi))\prec\phi$.
\end{description}
\end{propo}

\begin{demonstration}
Note that (HL1) is the exact analogue of (H1) in \cite{matu-kuhlm:hardy-deriv-gener-series},
replacing $\phi '/\phi$ by $l(\phi)$. The proof of the first statement is the exact analogue of the
 proof of \cite[Lemma 3.9]{matu-kuhlm:hardy-deriv-gener-series}, replacing $\phi '/\phi$ by $l(\phi)$.

Let $\alpha=\displaystyle\prod_{\phi\in\textrm{supp}\ \alpha}\phi^{\alpha_\phi} \in \Gamma$.
 Assume that $\alpha\succ 1$, thus
$\textrm{LE}(\alpha)>0$. By the first statement,
$l(\alpha)=\displaystyle\sum_{\phi\in\textrm{supp}\ \alpha}\alpha_\phi l(\phi)$.
Since by hypothesis $l$ is order preserving on $\Phi$, we have
$\textrm{LC}(l(\alpha))=\textrm{LE}(\alpha)\textrm{LC}(\phi_0)>0$ where $\phi_0=\textrm{LF}(\alpha)$, so $l(\alpha)>0$.

%The map $l$ obtained on $\Gamma$ is a pre-logarithmic section if and only if for any $\alpha\prec\beta$, $l(\alpha)<l(\beta)$.
%But, as is shown in \cite[Theorem 4.3]{matu-kuhlm:hardy-deriv-gener-series} for l'Hospital rule, this condition reduces to the same one on the fundamental monomials, which is (HL2).

For (HL3), we consider some monomial
$\alpha=\displaystyle\prod_{\phi\in\textrm{supp}\
\alpha}\phi^{\alpha_\phi}$ in $\Gamma\setminus\{1\}$. Then
$\textrm{LM}(l(\alpha))=\textrm{LM}(\displaystyle\sum_{\phi\in\textrm{supp}\
\alpha}\alpha_\phi l(\phi))=\textrm{LM}(\alpha_{\phi_0}l(\phi_0))$
where $\phi_0=\mbox{ LF }(\alpha)$ and $\alpha_{\phi_0}>0$. So
$l(\alpha)\prec\alpha$ for any $\alpha$ if and only if, for any
$\phi_0$ and $\alpha_{\phi_0}>0$,
$l(\phi_0)\prec\phi_0^{\alpha_{\phi_0}}$. This is equivalent to
(HL3).
\end{demonstration}

\begin{remark}
As in \cite[Corollaries 3.12 and 3.13]{matu-kuhlm:hardy-deriv-gener-series}, one can give analogous particular cases of (HL1).
\end{remark}

\subsection{Compatibility of pre-logarithms and derivations.}
We recall the following definition from \cite{matu-kuhlm:hardy-deriv-gener-series}:
\begin{defi}\label{defi:series-deriv}
A map $d:\mathds{K}\rightarrow\mathds{K}$, $a\mapsto a'$, verifying the following axioms is called a \textbf{series derivation}:
\begin{description}
\item[(D0)]  $1'=0$;
\item[(D1) Strong Leibniz rule:]$\forall\alpha=\displaystyle\prod_{\phi\in\textrm{supp}\
\alpha}\phi^{\alpha_\phi}\in\Gamma,\ (\alpha)'=
\alpha\displaystyle\sum_{\phi\in\textrm{supp}\ \alpha}\alpha_\phi\displaystyle\frac {\phi '}{\phi}$;
    \item[(D2) Strong linearity:]$\forall a=
\displaystyle\sum_{\alpha\in\textrm{Supp}\ a}a_\alpha
\alpha\in\mathds{K},\ a'= \displaystyle\sum_{\alpha\in\textrm{Supp}\
a}a_\alpha \alpha '$.
\end{description}
\end{defi}
Here we provide a criterion on the derivation to be compatible with the pre-logarithm:
\begin{defi}\label{defi:log-comp}
Let $(\mathds{K},l)$ be a pre-logarithmic field endowed with a
derivation $d$. Then $d$ is \textbf{log-compatible} if for all
$a\in\mathds{K}^{>0}$, we have $l(a)'=\displaystyle\frac{a'}{a}$. In
this case, we shall say the pre-logarithm $l$ is compatible with the
derivation $d$ or that $d$ and $l$ are compatible.
\end{defi}
In the case of a series morphism and a series derivation, it is sufficient to verify the
compatibility condition for the fundamental monomials:

\begin{propo}\label{propo:compat}
Let $(\mathds{K},l,d)$ be a generalized series field endowed with a series morphism $l$ and a series
 derivation $d$. Then  $d$ is log-compatible if and only if the following property holds:\begin{description}
    \item[(HL4)] $\forall \phi\in \Phi,\ l(\phi)'=\displaystyle\frac{\phi '}{\phi}$.
\end{description}
\end{propo}
\begin{demonstration}
Let $a=\alpha\ a_\alpha(1+\epsilon_a)\in\mathds{K}_{>0}$ where $\alpha=\displaystyle\prod_{\phi\in \textrm{supp }\alpha}\phi^{\alpha_\phi}$ and $a_\alpha\in\mathbb{R}_{>0}$. Using  (L), (D1), (D2), we compute: \\
$\begin{array}{rcl}
l(a)'&=&\left(l(\alpha)+\log(a_\alpha)+\displaystyle\sum_{k=1}^{+\infty}(-1)^{k-1} \epsilon_a^{k}\right)'\\
&=&\displaystyle\sum_{\phi\in\textrm{supp}\ \alpha}\alpha_\phi l(\phi)'+0+\left(\displaystyle\sum_{k=1}^{+\infty}(-1)^{k-1} \epsilon_a^{k-1}\right)\epsilon_a'
\end{array}$\mn
On the other hand, we compute: $a'=(\alpha\ a_\alpha (1+\epsilon_a))'=\alpha 'a_\alpha(1+\epsilon_a)+\alpha\ a_\alpha\epsilon_a'.$ Therefore:
$\begin{array}{lcl}
\displaystyle\frac{a'}{a}&=&\displaystyle\frac{\alpha '}{\alpha}+ \displaystyle\frac{\epsilon_a'}{1+\epsilon_a}
\end{array}$. Now, by (D1): $\alpha '=\alpha\displaystyle\sum_{\phi\in\textrm{supp}\ \alpha}\alpha_\phi\displaystyle\frac{\phi '}{\phi}$.\\ So : $\displaystyle\frac{a'}{a}=\displaystyle\sum_{\phi\in\textrm{supp}\ \alpha}\alpha_\phi\displaystyle\frac{\phi '}{\phi}+\displaystyle\frac{\epsilon_a'}{1+\epsilon_a}=\displaystyle\sum_{\phi\in\textrm{supp}\ \alpha}\alpha_\phi\displaystyle\frac{\phi '}{\phi}+ \left(\displaystyle\sum_{k=1}^{+\infty} (-1)^{k-1}\epsilon_a^{k-1}\right)\epsilon_a'$.\mn
Consequently: $l( a)'=\displaystyle\frac{a'}{a}$ if and only if $\displaystyle\frac{\phi '}{\phi}=l(\phi)'$ for all $\phi\in\textrm{supp}\ \alpha$.
\end{demonstration}

\section{Pre-logarithms and integration for Hardy type derivations.}\label{section:hardy-derivation}
We recall the following definition from \cite{matu-kuhlm:hardy-deriv-gener-series}:
\begin{defi}\label{defi:hardy_deriv} A derivation
$d$ on $\mathds{K}$ is a {\bf Hardy type derivation} if :
\begin{description}
\item[(HD1)] the \textbf{sub-field of constants} of
$\mathds{K}$ is $\mathbb{R}$.
\item[(HD2)] $d$ verifies \textbf{l'Hospital's rule} :
$\forall a,b\in \mathds{K}^*$ with $a,b \nasymp 1$ we have\n
$a\preccurlyeq b\Leftrightarrow a'\preccurlyeq b'$.
\item[(HD3)] the logarithmic derivation is \textbf{compatible with the dominance relation}:
 $\forall  a,b\in\mathds{K}$ with $|a|\succ|b|\succ 1$,
we have $\displaystyle\frac{a'}{a}\succcurlyeq
\displaystyle\frac{b'}{b}$. Moreover,
$\displaystyle\frac{a'}{a}\asymp\displaystyle\frac{b'}{b}$ if and
only if $a$ and $b$ are comparable.
\end{description}
\end{defi}
%-----------------------------------------------------------------------------

\subsection{The monomial asymptotic integral.}\label{section:asymp-integ}

\textit{For the rest of this section, we assume that $d$ is a Hardy type series  derivation on $\mathds{K}$.}
\sn
\begin{notation}
 Set \[\forall\phi\in\Phi,\> \theta^{(\phi)}:=\mbox{LM}\left(\displaystyle\frac{\phi
'}{\phi}\right)\>,\>\>\Theta: =\left\{\theta^{(\phi)},\
\phi\in\Phi\right\}\>,\>\>
 \mbox{ and }\>\>
 \hat{\theta}\>:=\> \mbox{ g.l.b. }_{\preccurlyeq}\Theta\]
  if it exists in $\Gamma$.\sn
  Adopting the notation of \cite{rosenlicht:rank}, we write below:
$\Psi:=\left\{\mbox{LM}\left(\displaystyle\frac{a'}{a}\right);\ a\in
\mathds{K}^*,\ a\nasymp\ 1\right\}$.
\end{notation}
We will make use of the following result \cite[Theorem 4.3;
Corollary 4.4]{matu-kuhlm:hardy-deriv-gener-series}: \mn \textit{A
series derivation on $\mathds{K}$ is of Hardy type if and only if
the following condition holds:\begin{description}
\item[(H3')] $\forall\phi\prec\psi\in\Phi$,
$\theta^{(\phi)}\prec\theta^{(\psi)}$ and $\mbox{ LF
}\left(\displaystyle\frac{\theta^{(\phi)}}{\theta^{(\psi)}}\right)\prec\psi$.
\end{description}}

\begin{defi}\label{defi:asymp-integ}
 We say that $b\in\mathds{K}$ is an \textbf{asymptotic integral} of $a\in\mathds{K}$ if $b'\sim a$, equivalently if $ b'\sim \mbox{LT}(a)$.
We say that $b$ is an \textbf{integral} of $a$ if $b'=a$.
\end{defi}

\begin{theo}\label{theo:integ-asymp-prelog}
 A series $a\in\mathds{K}^*$ has an asymptotic
integral if and only if  $a\nasymp\ \mbox{ g.l.b.
}_\preccurlyeq\Psi$.
\end{theo}
This theorem is proved for Hardy fields in \cite[Proposition 2 and Theorem 1]{rosenlicht:rank}.
 As noted in \cite{matu-kuhlm:hardy-deriv-gener-series}, it suffices to observe that Rosenlicht's proof
 only uses the properties of what we call a Hardy type derivation in Definition \ref{defi:hardy_deriv}.
 If $d$ is moreover a series derivation, it suffices to consider fundamental monomials as we establish below.

\begin{propo}\label{propo:theta-hat}
Assume that $d$ is a Hardy type series derivation on $\mathds{K}$.
Let $a\in\mathds{K}^*$ with $a\nasymp 1$. Then $$\mathrm{LT}\left(\displaystyle\frac{a '}{a}\right)=
\mathrm{LE}(a)\mathrm{LT} \left(\displaystyle\frac{\mathrm{LF}(a)'}{\rm{LF}(a)}\right).$$
More precisely $\mathrm{LM}\left(\displaystyle\frac{a '}{a}\right)=
\theta^{(\rm{LF}(a))}$ and $\mathrm{LC}\left(\displaystyle\frac{a '}{a}\right) =\mathrm{LE}(a)\mathrm{LC} \left(\displaystyle\frac{\mathrm{LF}(a)'}{\rm{LF}(a)}\right)$. \sn In particular, $\hat{\theta}= \mbox{ g.l.b.
}_\preccurlyeq\Psi$.
\end{propo}
\begin{demonstration}
Let $1\nasymp a=a_\alpha\alpha+\cdots\in\mathds{K}^*$ with
 $\alpha=\displaystyle\prod_{\phi\in\mathrm{supp}\ \alpha}\phi^{\alpha_\phi}$.
Set $\phi_0=\mbox{LF}(a)=\mbox{LF}(\alpha)$ and $\alpha_{\phi_0}= \mathrm{LE}(a)=\mathrm{LE}(\alpha)$.
 We compute: $$a'=a_\alpha\alpha '+\cdots=
 a_\alpha\alpha\left(\alpha_{\phi_0}\displaystyle\frac{\phi_0 '}{\phi_0}+\cdots\right)+\cdots=
 (a_\alpha\alpha_{\phi_0})\alpha\displaystyle\frac{\phi_0 '}{\phi_0}+\cdots$$ Therefore:
  %$\mbox{LM}(a')=\mbox{LM}\left(\alpha\displaystyle\frac{\phi_0 '}{\phi_0}\right)=\alpha\theta^{(\phi_0)}$.
$$\mathrm{LT}\left(\displaystyle\frac{a '}{a}\right)= \displaystyle\frac{\mathrm{LT}(a')}{ \mathrm{LT}(a)}=
\displaystyle\frac{(a_\alpha\alpha_{\phi_0}) \alpha\mathrm{LT} \left(\displaystyle\frac{\phi_0 '}{\phi_0}\right)
}{a_\alpha\alpha}= \alpha_{\phi_0}\mathrm{LT} \left(\displaystyle\frac{\phi_0 '}{\phi_0}\right).$$.
\end{demonstration}

\cite[Theorem 1]{rosenlicht:rank} gives a parametrized family of
asymptotic integrals of an (asymptotically integrable) element $a$.
For a Hardy type series derivations, we compute in Proposition
\ref{propo:exist-integ-asymp} below a specific asymptotic integral,
which turns out to be a non monic monomial (i. e. of the form
$r\alpha$ with $r\in\mathbb{R}$ and $\alpha \in \Gamma$), uniquely
determined by $a$.

\begin{notation}
 We call the asymptotic integral computed in Proposition \ref{propo:exist-integ-asymp} below the {\bf monomial asymptotic integral} of $a$, and denote it by $\mathrm{a.i.}(a)$.
\end{notation}

\begin{lemma}\label{lemme:psi}
Let $\alpha\in\Gamma$ with $\alpha\not= \hat{\theta}$. There exists
a uniquely determined fundamental monomial $\psi_\alpha\in\Phi$
which satisfies
$\mathrm{LF}\left(\displaystyle\frac{\alpha}{\theta^{(\psi_\alpha)}}\right)=
\psi_\alpha$.
\end{lemma}
\begin{demonstration}
First, suppose that  $\alpha\succ \hat{\theta}$. Take a monomial
$\beta\succ 1$ with $\alpha\succ \displaystyle\frac{\beta
'}{\beta}$. Set  $\phi:=\mathrm{LF}(\beta)$, so $\displaystyle\frac{\beta
'}{\beta}\asymp \theta^{(\phi)}$ by Proposition \ref{propo:theta-hat}.
Set $\beta_0:=\min\left\{\beta,
\displaystyle\frac{\alpha}{\theta^{(\phi)}}\right\}$ and
$\phi_0:=\mathrm{LF}(\beta_0)$. Since $\beta\succcurlyeq\beta_0\succ
1$, we have $\phi\succcurlyeq\phi_0$, so
$\theta^{(\phi)}\succcurlyeq\theta^{(\phi_0)}$. We deduce that
$\alpha\succ\theta^{(\phi_0)}$ and
$\displaystyle\frac{\alpha}{\theta^{(\phi_0)}}\succcurlyeq
\displaystyle\frac{\alpha}{\theta^{(\phi)}}\succcurlyeq \beta_0\succ
1$. If we set
$\phi_1:=\mathrm{LF}\left(\displaystyle\frac{\alpha}{\theta^{(\phi_0)}}\right)$,
then $\phi_1\succcurlyeq\phi_0$. We compute:
$\mathrm{LF}\left(\displaystyle\frac{\alpha}{\theta^{(\phi_1)}}\right)
= \mathrm{LF}\left(\displaystyle\frac{\alpha}{\theta^{(\phi_0)}}.
\displaystyle\frac{\theta^{(\phi_0)}}{\theta^{(\phi_1)}}\right)$. By (H3'):
$\mathrm{LF}\left(\displaystyle\frac{\theta^{(\phi_0)}}{\theta^{(\phi_1)}}\right)
\prec \phi_1$. We obtain:
$\mathrm{LF}\left(\displaystyle\frac{\alpha}{\theta^{(\phi_1)}}\right)
=\max\left\{\mathrm{LF}\left(\displaystyle\frac{\alpha}{\theta^{(\phi_0)}}\right); \mathrm{LF}\left(\displaystyle\frac{\theta^{(\phi_0)}}{\theta^{(\phi_1)}}\right)\right\} =\phi_1$.
Set $\psi_\alpha:=\phi_1$.
\sn Now suppose that   $\alpha\prec
\hat{\theta}$. Let $\alpha_1\in\Gamma$ such that $\alpha\prec\alpha_1\preccurlyeq
\hat{\theta}$. Set $\phi_0:=\mathrm{LF}\left(\displaystyle\frac{\alpha}{\alpha_1}\right)$, then
$\displaystyle\frac{\alpha}{\theta^{(\phi_0)}}=
\displaystyle\frac{\alpha}{\alpha_1}.
\displaystyle\frac{\alpha_1}{\theta^{(\phi_0)}}\preccurlyeq
\displaystyle\frac{\alpha}{\alpha_1}\prec 1$. Set
$\phi_1:=\mathrm{LF}\left(\displaystyle\frac{\alpha}{\theta^{(\phi_0)}}\right)$.
We deduce that $\phi_1\succcurlyeq\phi_0$, and compute
$\mathrm{LF}\left(\displaystyle\frac{\alpha}{\theta^{(\phi_1)}}\right)=\phi_1$ as above.
Set $\psi_\alpha:=\phi_1$. This concludes the proof of the existence of $\psi_\alpha$.
\sn Consider now a monomial $\alpha\nasymp \hat{\theta}$, and denote $\psi_1$
and $\psi_2$ two fundamental monomials such that
$\mbox{LF}\left(\displaystyle\frac{\alpha}{\theta^{(\psi_i)}}\right)=\psi_i$ for $i=1,2$.
Assume for instance that $\psi_1\prec\psi_2$. We would have
$\mbox{LF}\left(\displaystyle\frac{\alpha}{\theta^{(\psi_2)}}\right)=
\mbox{LF}\left(\displaystyle\frac{\alpha}{\theta^{(\psi_1)}}.
\displaystyle\frac{\theta^{(\psi_1)}}{\theta^{(\psi_2)}}\right)=
\psi_2$. Since $\mbox{LF}\left(\displaystyle\frac{\alpha}{\theta^{(\psi_1)}}\right)=\psi_1$, we would have
$\mbox{LF}\left(\displaystyle\frac{\theta^{(\psi_1)}}{\theta^{(\psi_2)}}\right)
=\psi_2$, which contradicts (H3').
\end{demonstration}

\begin{propo}\label{propo:exist-integ-asymp}
Let $a\in\mathds{K}^*$ with $a\nasymp \hat{\theta}$, and set $\alpha:=\mathrm{LM}(a)$. Then:
\begin{center}
$\mathrm{a.i.}(\alpha)= \displaystyle\frac{\alpha}{\mathrm{LE}\left(\displaystyle
\frac{ \alpha}{\theta^{(\psi_\alpha)}}\right)
\mathrm{LT}\left(\displaystyle\frac{\psi_\alpha '}{\psi_\alpha}\right)}\ \ \ \ \  $ and
$\ \ \ \ \ \mathrm{a.i.}(a)=\mathrm{LC}(a) \mathrm{a.i.}(\alpha)$
\end{center}
\end{propo}
\begin{demonstration} %As was already noticed after Definition \ref{defi:asymp-integ}, it suffices to consider a monomial $\alpha\nasymp \hat{\theta}$.
%, and assume that there is no $\psi\in\Phi$ such that
% $\mbox{LF}\left(\displaystyle\frac{\alpha}{\theta^{(\psi)}}\right)=\psi$.
Below, set $m:=\mathrm{a.i.}(\alpha)=
\displaystyle\frac{\alpha}{\mathrm{LE}\left(\displaystyle \frac{
\alpha}{\theta^{(\psi_\alpha)}}\right)
\mathrm{LC}\left(\displaystyle\frac{\psi_\alpha
'}{\psi_\alpha}\right)\theta^{(\psi_\alpha)}}$. \sn Since
$\mathrm{LF}(m)=\mathrm{LF}\left(\displaystyle\frac{\alpha}{\theta^{(\psi_\alpha)}}\right)
=\psi_\alpha$, using Proposition \ref{propo:theta-hat}, we compute:
$$ \mathrm{LT}\left(\displaystyle\frac{m'}{m}\right)=\mathrm{LE}(m) \mathrm{LT}\left(\displaystyle\frac{\psi_\alpha '}{\psi_\alpha}\right)$$
Since
$\mathrm{LE}(m)=\mathrm{LE}\left(\displaystyle\frac{\alpha}{\theta^{(\psi_\alpha)}}\right)\>,$
we compute: \[\mathrm{LT}(m')=m\mathrm{LE}(m)
\mathrm{LT}\left(\displaystyle\frac{\psi_\alpha
'}{\psi_\alpha}\right)
=\displaystyle\frac{\alpha}{\mathrm{LE}\left(\displaystyle\frac{
\alpha}{\theta^{(\psi_\alpha)}}\right)\mathrm{LT}\left(\displaystyle
\frac{\psi_\alpha
'}{\psi_\alpha}\right)}.\mathrm{LE}\left(\displaystyle\frac{
\alpha}{\theta^{(\psi_\alpha)}}\right)\mathrm{LT}\left(\displaystyle
\frac{\psi_\alpha '}{\psi_\alpha}\right) =\alpha\>,\] as desired.
\sn Denote $b:=\mathrm{a.i.}(a)$. We have:
$\mathrm{LT}(b')=\mathrm{LT}\left(\mathrm{LC}(a)m'\right)=\mathrm{LC}(a)\mathrm{LT}\left(m'\right)=
\mathrm{LC}(a)\alpha=\mathrm{LT}(a)$, as desired.
\end{demonstration}

\begin{notation}
In the sequel, to simplify the notations, we will write $\psi$ instead of $\psi_\alpha$ (of Lemma \ref{lemme:psi}) if the context is clear.
\end{notation}

\subsection{Constructing pre-logarithms as anti-derivatives.}
 In the following theorem, we give a criterion for
$(\mathds{K},d)$ to carry a pre-logarithm,  compatible with the
derivation. Moreover, we will require this pre-logarithm to be
induced by a pre-logarithmic section which is a series morphism. The
construction relies on the computation of the anti-derivatives of
$\displaystyle\frac{\phi '}{\phi}$, $\phi\in\Phi$.
\begin{theo}\label{theo:existence-prelog}
Let $d$ be a Hardy type series derivation on $\mathds{K}$.
There exists a unique pre-logarithmic section $l$ on
$\mathds{K}$ which is a series morphism, for which the induced
pre-logarithm is compatible with the derivation, if and only if the following two conditions hold:\begin{enumerate}
    \item $\hat{\theta}\notin\displaystyle\bigcup_{\phi\in\Phi}\textrm{Supp}\ \displaystyle\frac{\phi '}{\phi}$;
\item $\forall\phi\in\Phi$, $\forall\tau^{(\phi)}\in\textrm{Supp}\ \displaystyle\frac{\phi '}{\phi}$,  $\mathrm{a.i.}(\tau^{(\phi)})\succ 1$.
% $\displaystyle\frac{\tau^{(\phi)}}{\theta^{(\psi)}}\succ 1$ where $\psi:=\psi_{\tau^{(\phi)}}$ (Lemma \ref{lemme:psi}). In particular, $\psi\prec\phi$.
\end{enumerate}
 Moreover, this pre-logarithm verifies (GA).
\end{theo}
\begin{demonstration} To define a pre-logarithm $l$ on $\mathds{K}_{>0}$, it suffices to define a pre-logarithmic section $l$ on $\Gamma$. We set $l(1):=0$. By (D1), for
any $\alpha=\displaystyle\prod_{\phi\in\textrm{supp}\
\alpha}\phi^{\alpha_\phi}\in\Gamma\backslash\{1\}$, we have
$\displaystyle\frac{\alpha '}{\alpha}=
\displaystyle\sum_{\phi\in\textrm{supp}\
\alpha}\alpha_\phi\displaystyle\frac{\phi '}{\phi}$. Assume that for
any $\phi\in\Phi$, there exists $l(\phi)\in\mathds{K}^{\succ 1}$
such that (HL4) holds, i.e. $l(\phi)'=\displaystyle\frac{\phi '}{\phi}$. (The proof of the existence of such $l(\phi)\in\mathds{K}^{\succ 1}$ will be
established below).
We apply Proposition \ref{propo:pre-log} to extend $l$ to a series morphism on $\Gamma$.
Suppose, as in (HL1), that there
exists a strictly decreasing sequence
$(\phi_n)_{n\in\mathbb{N}}\subset\Phi$ and an increasing sequence
 $(\lambda^{(n)})_{n\in\mathbb{N}}\subset\Gamma$ such that for any $n$, $\lambda^{(n)}\in\textrm{Supp}\ l(\phi_n)$.
By (HD2), $\tau^{(n)}:=\mathrm{LM}((\lambda^{(n)})')$ defines an increasing sequence in $\Gamma$ such that for any $n$,
 $\tau^{(n)}\in\textrm{Supp}\ \displaystyle\frac{\phi_n '}{\phi_n}$.
  This implies that \cite[ (H1)]{matu-kuhlm:hardy-deriv-gener-series} holds,
   contradicting the fact that $d$ is a series derivation.
   Therefore, for any $\alpha\in\Gamma$, we can indeed define $l(\alpha):=\displaystyle\sum_{\phi\in\textrm{supp}\ \alpha}\alpha_\phi l(\phi)$.
\sn   Note that by (HD2), (HL2) holds. Thus $l$ would be the pre-logarithmic section induced  by the given $l$ on $\Phi$. Furthermore, this series morphism $l$ is compatible with the derivation (Proposition \ref{propo:compat}).
\mn It remains to prove the existence of such
$l(\phi)\in\mathds{K}^{\succ 1}$.
We adapt to our context \cite[Theorem 1]{kuhl:Hensel-lemma}, with the "spherically complete" ultrametric space $(\mathds{K},u)$ where $u(a,b):=\mathrm{LM}(a-b)$, and the map $f:=d$.
\begin{lemma}[\cite{kuhl:Hensel-lemma}, Theorem 1]
Let $\phi\in\Phi$. We suppose that for any $a\in\mathds{K}$ with $a'\neq \displaystyle\frac{\phi '}{\phi}$, there exists $b\in\mathds{K}$ such that:
\begin{description}
    \item[(AT1)] $\mathrm{LM}\left(b'-\displaystyle\frac{\phi '}{\phi}\right)\succ \mathrm{LM}\left(a'-\displaystyle\frac{\phi '}{\phi}\right)$;
\item[(AT2)] $\forall c\in\mathds{K}$, if $\mathrm{LM}\left(a-c\right)\succ \mathrm{LM}\left(a-b\right)$, then $\mathrm{LM}\left(a'-c'\right)\succ \mathrm{LM}\left(a'-\displaystyle\frac{\phi '}{\phi}\right)$.
\end{description}
Then there exists $l(\phi)\in\mathds{K}$ such that $l(\phi)'=\displaystyle\frac{\phi '}{\phi}$.
\end{lemma}
\begin{demonstration}
Let $a\in\mathds{K}$. By (D1) and (D2), we can denote $\mathrm{LT}\left(a'-\displaystyle\frac{\phi '}{\phi}\right)=c_0\alpha\tau^{(\tilde{\phi})}$ for some $c_0\in\mathbb{R}$, $\alpha\in\mathrm{Supp}\ a\cup \{1\}$, $\tilde{\phi}\in\mathrm{supp}\ \alpha\cup\{\phi\}$ and $\tau^{(\tilde{\phi})}\in\mathrm{Supp}\ \displaystyle\frac{\tilde{\phi} '}{\tilde{\phi}}$.
\begin{claim}\label{claim:a.i.-derive}
Provided Hypothesis 1 and 2 of Theorem \ref{theo:existence-prelog}, we consider $\alpha\in\Gamma$, $\alpha\nasymp 1$. Then any monomial $\beta=\alpha\tau^{(\tilde{\phi})}\in \mathrm{Supp}\ (\alpha ')$ (where $\tilde{\phi}\in\mathrm{supp}\ \alpha$ and $\tau^{(\tilde{\phi})}\in\mathrm{Supp}\ \displaystyle\frac{\tilde{\phi} '}{\tilde{\phi}}$ by (D1)) admits an asymptotic integral. Moreover, $\psi_\beta=\mathrm{LF}(\alpha)$ and $\mathrm{LE}(\beta)=\mathrm{LE}(\alpha)$.
\end{claim}
Indeed, by Lemma \ref{lemme:psi} and Proposition \ref{propo:exist-integ-asymp}, we show that $\mathrm{LF}\left(\displaystyle\frac{\alpha\tau^{(\tilde{\phi})} }{\theta^{(\psi)}}\right) =\psi$. Set $\psi:=\mathrm{LF}(\alpha)$, therefore $\psi\succcurlyeq\tilde{\phi}$. Denote by $\tilde{\psi}$ the unique fundamental monomial such that $\mathrm{LF}\left(\displaystyle\frac{\tau^{(\tilde{\phi})} }{\theta^{(\tilde{\psi})}} \right)=\tilde{\psi}$ (which exists since $\tau^{(\tilde{\phi})}\nasymp\hat{\theta}$ by Hypothesis 1). Since $\displaystyle\frac{\tau^{(\tilde{\phi})} }{\theta^{(\tilde{\psi})}}\succ 1$ by Hypothesis 2, we have $\tilde{\psi}\prec\tilde{\phi}$. Consequently, $\tilde{\psi}\prec\psi$, so $\mathrm{LF}\left(\displaystyle\frac{\theta^{(\tilde{\psi})} }{\theta^{(\psi)}}\right) \prec \psi$ by (H3'). Then, using the ultrametric triangular inequality for LF, we compute:
$$\mathrm{LF}\left(\displaystyle\frac{\alpha\tau^{(\tilde{\phi})} }{\theta^{(\psi)}}\right)=\mathrm{LF}\left(\alpha\displaystyle\frac{\tau^{(\tilde{\phi})} }{\theta^{(\tilde{\psi})}}\displaystyle\frac{\theta^{(\tilde{\psi})} }{\theta^{(\psi)}}\right)=\mathrm{LF}(\alpha)=\psi\ \ \ \ \  \mathrm{and}\ \ \ \ \mathrm{LE}(\beta)=\mathrm{LE}(\alpha).$$

\noindent Consequently, $c_0\alpha\tau^{(\tilde{\phi})}$ admits an asymptotic integral monomial. To conclude the proof of (AT1), it suffices to set $b:=a-\mathrm{a.i.}(c_0\alpha\tau^{(\tilde{\phi})})$.
\sn
Concerning (AT2), we consider $c\in\mathds{K}$ such that $$\mathrm{LM}\left(a-c\right)\succ \mathrm{LM}\left(a-b\right)=\mathrm{LM}\left(\mathrm{a.i.} (c_0\alpha\tau^{(\tilde{\phi})})\right)= \displaystyle\frac{\alpha\tau^{(\tilde{\phi})} }{\theta^{(\psi)}}.$$
By (HD2), we have:
$$\mathrm{LM}\left(a'-c'\right)\succ \mathrm{LM}\left[\left(\displaystyle\frac{\alpha\tau^{(\tilde{\phi})} }{\theta^{(\psi)}}\right)'\right] =\alpha\tau^{(\tilde{\phi})}= \mathrm{LM}\left(a'-\displaystyle\frac{\phi '}{\phi}\right).$$
\end{demonstration}
Note that $l(\phi)$ is defined up to addition by a real constant. \textit{We choose the $l(phi)$'s so that this real constant is zero, i.e. $1\notin\mathrm{Supp}\ l(\phi)$}.
\sn
We prove now that $l(\phi)\in\mathds{K}^{\succ 1}$ for any $\phi\in\Phi$. Suppose not, and denote by $\lambda^{(\phi)}$ the greatest monomial in $\mathrm{Supp}\ l(\phi)$ such that $\lambda\prec 1$. Then, $\mathrm{LM}\left(\lambda '\right)=\lambda\theta^{(\psi)}$, where $\psi=\mathrm{LF}(\lambda)$. We consider two cases. Either $\lambda\theta^{(\psi)}=\tau\in\mathrm{Supp}\ \displaystyle\frac{\phi '}{\phi}$, which is impossible since $\mathrm{a.i.}(\tau)\succ 1$ by Hypothesis 2. Or $\lambda\theta^{(\psi)}=\tilde{\lambda}\tilde{\tau}$ for some $\tilde{\lambda}\succ 1$, $\tilde{\phi}\in\mathrm{supp}\ \tilde{\lambda}$ and $\tilde{\tau}\in\mathrm{Supp}\ \displaystyle\frac{\tilde{\phi} '}{\tilde{\phi}}$, meaning that, up to multiplication by a real coefficient, $\lambda$ is the asymptotic integral monomial of $\tilde{\lambda}\tilde{\tau}$. But, computing  $\mathrm{a.i.}\left(\tilde{\lambda}\tilde{\tau}\right)$ as in the proof of (AT1) in the preceding lemma, we obtain:
$$ \mathrm{LM}\left[\mathrm{a.i.}\left(\tilde{\lambda}\tilde{\tau}\right)\right] =\displaystyle\frac{\tilde{\lambda}\tau^{(\tilde{\phi})} }{\theta^{(\psi)}}$$ with $\psi:=\mathrm{LF}\left(\tilde{\lambda}\right)=\mathrm{LF}\left( \displaystyle\frac{\tilde{\lambda}\tau^{(\tilde{\phi})} }{\theta^{(\psi)}}\right)$ and $\mathrm{LE}\left( \displaystyle\frac{\tilde{\lambda}\tau^{(\tilde{\phi})} }{\theta^{(\psi)}}\right)=\mathrm{LE}\left(\tilde{\lambda}\right)>0$. This means that $\mathrm{a.i.}\left(\tilde{\lambda}\tilde{\tau}\right)\succ 1$: contradiction.
\mn
To conclude the proof of the theorem, we show that the
pre-logarithm is uniquely determined, and that it verifies (GA).
Indeed, let  $l_1$ and $l_2$ be two pre-logarithms compatible with
$d$, and $a\in\mathds{K}_{>0}$. So
$l_1(a)'=\displaystyle\frac{a'}{a}=l_2(a)'$, which means that
$l_1(a)=l_2(a)+c$ for some $c\in\mathbb{R}$. But if we take $a=1$,
then $l_1=l_2=\log$, so $l_1(1)=l_2(1)=0$ which implies that $c=0$.

Concerning (GA), since the derivation verifies l'Hospital's rule
(HD2), we observe that, for any $\phi$, the leading monomial of
$l(\phi)$ is $\displaystyle\frac{\theta^{(\phi)}}{\theta^{(\psi)}}$
where $\psi$ is the fundamental monomial such that
$\textrm{LF}\left(\displaystyle\frac{\theta^{(\phi)}}{\theta^{(\psi)}}\right)
=\psi$ (exists by Hypothesis 1.). By Hypothesis 2., we have moreover that
$\displaystyle\frac{\theta^{(\phi)}}{\theta^{(\psi)}}\succ
1$, so $\phi\succ\psi$.
Thus, we obtain $\textrm{LF}(l(\phi))=\psi\prec \phi$, as
desired.
\end{demonstration}
In the next result and in its proof, we give a description of the $l(\phi)$'s via a method that we may call an \textbf{iterated asymptotic integration}.

\begin{coro}\label{next}
With the same hypothesis as in Theorem \ref{theo:existence-prelog}, for any $\phi\in\Phi$, if we denote
$l(\phi)=\displaystyle\sum_{\lambda\in\textrm{Supp}\ l(\phi)}d_{\lambda}\lambda\in\mathds{K}^{\succ 1}$, then for any $\lambda\in\textrm{Supp}\ l(\phi)$, there is $n\in\mathbb{N}$ such that:
\begin{center}
$\lambda=\displaystyle\prod_{i=1}^n \displaystyle\frac{\tau^{(\phi_{i})}}{\theta^{(\psi)}}\ \ \ \ \ $ and  $\ \ \ \ \  d_\lambda =\displaystyle\frac{\prod_{i=1}^n c_{\tau^{(\phi_{i})}}}{ \left(\beta_0c_{0,\psi}\right)^n} $
\end{center}
where:
\begin{description}
\item[a)] $\tau^{(\phi_{1})}=\tau^{(\phi)}\in  \mathrm{Supp}\ \displaystyle\frac{\phi '}{\phi}$ and $\psi=\psi_{\tau^{(\phi)}}$ (Lemma \ref{lemme:psi}: i.e. $\psi$ verifies $\mathrm{LF}\left(\displaystyle\frac{\tau^{(\phi)}}{\theta^{(\psi)}}\right)=\psi$);
\item[b)] for any $i=2,\ldots,n$, $c_{\tau^{(\phi_{i})}}\tau^{(\phi_i)}$ is a monomial of $\displaystyle\frac{\phi_i '}{\phi_i}$ for some $\phi_i\preccurlyeq\psi$ with $\tau^{(\phi_i)}\prec\theta^{(\psi)}$;
\item[c)]  $\beta_0=\mathrm{LE}\left(\displaystyle\frac{ \tau^{(\phi)}}{\theta^{(\psi)}}\right)>0$ and $c_{0,\psi}=\mathrm{LC} \left(\displaystyle\frac{\psi '}{\psi}\right)$;
\item[d)] for any $k=1,\ldots,n$, $\mathrm{LF}\left(\displaystyle\prod_{i=1}^k \displaystyle\frac{\tau^{(\phi_{i})}}{\theta^{(\psi)}}\right)=\psi$
and $\mathrm{LE}\left(\displaystyle\prod_{i=1}^k \displaystyle\frac{\tau^{(\phi_{i})}}{\theta^{(\psi)}}\right)
=\beta_0 >0$.
\end{description}
\end{coro}

%{\it To simplify the notations, from now until the end of the proof, we ``forget" the coefficients of the monomials considered in the computations (we just need their existence, which will appear clearly through the asymptotic integration process)}. Assumptions 1. and 2. of the theorem will be used in the following:

\begin{demonstration} Let $\phi\in\Phi$. We set the iterated asymptotic integration of $\displaystyle\frac{\phi '}{\phi}$ as being the fixed point of the following map $f$ (we prove below that such a fixed point is well defined, unique and equal to $l(\phi)$). Given a series $l=\displaystyle\sum_{\lambda\in S } d_{\lambda}\lambda$  (which can be thought as an approximation of $l(\phi)$), by (D1) and (D2) we have:
$$\left[l(\phi)-l\right]'=\displaystyle\frac{\phi '}{\phi}- l'=\displaystyle\sum_{\lambda\in\textrm{Supp}\ l(\phi)\setminus S } \displaystyle\sum_{\tilde{\phi}\in\textrm{supp}\ \lambda } \displaystyle\sum_{\tau^{(\tilde{\phi})}\in\textrm{Supp}\ \phi '/\phi } (d_\lambda \tilde{c}_{\tilde{\phi}}).\lambda\tau^{(\tilde{\phi})}.$$
Since any of the terms $(d_\lambda \tilde{c}_{\tilde{\phi}}).\lambda\tau^{(\tilde{\phi})}$ admits an asymptotic integral monomial (Claim \ref{claim:a.i.-derive}), we set:
$$\mathrm{A.I.}\left(\left[l(\phi)-l\right]'\right):=\displaystyle\sum_{\lambda\in\textrm{Supp}\ l(\phi)\setminus S } \displaystyle\sum_{\tilde{\phi}\in\textrm{supp}\ \lambda } \displaystyle\sum_{\tau^{(\tilde{\phi})}\in\textrm{Supp}\ \phi '/\phi } \mathrm{a.i.}\left[(d_\lambda \tilde{c}_{\tilde{\phi}}).\lambda\tau^{(\tilde{\phi})}\right]\ \ \ \mathrm{and} \ \ \ \mathrm{A.I.}(0):=0.$$
and
$$f(l) :=l + \mathrm{A.I.}\left(\left[l(\phi)-l\right]'\right) $$
Note that $l(\phi)$ is a fixed point for $f$. We adapt to our context \cite[Theorem 1]{p-crampe-rib_FP-gen-series} for the ultrametric $u(a,b):=\mathrm{LM}(a-b)$), provided the fact that $(\mathds{K},u)$ is spherically complete:
\begin{lemma}[\cite{p-crampe-rib_FP-gen-series}, Theorem 1]
Since $\mathds{K}$ is spherically complete, any contracting map $f:\mathds{K}\rightarrow\mathds{K}$ has exactly one fixed point.
\end{lemma}
Our map $f$ is contracting. Indeed, given $l_1,l_2\in \mathds{K}$, $l_1\neq l_2$, we compute:
$$\begin{array}{lcl}
f(l_1)-f(l_2)&=&l_1 - \mathrm{A.I.}\left(\displaystyle\frac{\phi '}{\phi}-l_1'\right) -l_2+ \mathrm{A.I.}\left(\displaystyle\frac{\phi '}{\phi}-l_2'\right)\\
&=&(l_1-l_2)-\mathrm{A.I.}\left[(l_1-l_2)'\right].
\end{array}$$
Therefore: $u[f(l_1)-f(l_2)]=\mathrm{LM}\left[(l_1-l_2)-\mathrm{A.I.}\left[(l_1-l_2)'\right] \right] <\mathrm{LM}(l_1-l_2)=u(l_1,l_2)$.
Consequently, $l(\phi)$ is the unique fixed point of $f$.
\mn
To obtain the desired properties for $l(\phi)$, we proceed by induction along the iterated asymptotic integration. We begin with $l=0$. Thus, we compute the asymptotic integral of any monomial $c_{\tau^{(\phi)}}\tau^{(\phi)}$ of $\displaystyle\frac{\phi '}{\phi}$.  By Proposition \ref{propo:exist-integ-asymp} and Hypothesis 1, its asymptotic integral exists and is of the form: $$d_\lambda\lambda:=\displaystyle\frac{c_{\tau^{(\phi)}}}{ \beta_0c_{0,\psi}} \displaystyle\frac{\tau^{(\phi)}}{\theta^{(\psi)}}$$
 where $\psi:=\psi_{\tau^{(\phi)}}$, $\beta_0:=\mathrm{LE}\left(\displaystyle\frac{ \alpha}{\theta^{(\psi)}}\right)$ and $c_{0,\psi}:=\mathrm{LC} \left(\displaystyle\frac{\psi '}{\psi}\right)$.
Moreover, by Hypothesis 2, $\lambda\succ 1$ as desired.
\sn
We consider now  $f^n(0)$ for some $n\in\mathbb{N}$ which we denote by the series $l=\displaystyle\sum d_\lambda\lambda$, supposing that properties a) to d) hold for it. Then any term in $\left[l(\phi)-l\right]'$ is of the form $$(d_\lambda \tilde{c}_{\tilde{\phi}}).\lambda\tau^{(\tilde{\phi})}=\displaystyle\frac{\left( \prod_{i=1}^n c_{\tau^{(\phi_{i})}}\right)\tilde{c}_{\tilde{\phi}}}{ \left(\beta_0c_{0,\psi}\right)^n} \left(\displaystyle\prod_{i=1}^n \displaystyle\frac{\tau^{(\phi_{i})}}{\theta^{(\psi)}}\right)\tau^{(\tilde{\phi})}$$ where $\tilde{\phi}\in\mathrm{supp}\ \lambda$, so $\tilde{\phi}\preccurlyeq\psi$,   and $\tilde{c}_{\tilde{\phi}}\tau^{(\tilde{\phi})}$ is a monomial of $\mathrm{Supp}\ \displaystyle\frac{\tilde{\phi} '}{\tilde{\phi}}$ with $\tau^{(\tilde{\phi})}\prec \theta^{(\psi)}$. By Proposition \ref{propo:exist-integ-asymp}, Claim \ref{claim:a.i.-derive} and the induction hypothesis, its asymptotic integral is:
$$ d_{\tilde{\lambda}}\tilde{\lambda}:=\displaystyle\frac{d_\lambda \tilde{c}_{\tilde{\phi}}}{ \beta_0c_{0,\psi}} \displaystyle\frac{\lambda}{\theta^{(\psi)}} =\displaystyle\frac{\left( \prod_{i=1}^{n+1} c_{\tau^{(\phi_{i})}}\right)}{ \left(\beta_0c_{0,\psi}\right)^{n+1}} \displaystyle\prod_{i=1}^{n+1} \displaystyle\frac{\tau^{(\phi_{i})}}{\theta^{(\psi)}}$$
where $\phi_{n+1}:=\tilde{\phi}$, $\displaystyle\tau^{(\phi_{n+1})}:=\tau^{(\tilde{\phi})}$ and $c_{\tau^{(\phi_{n+1})}}:=\tilde{c}_{\tilde{\phi}}$. Note that $\mathrm{LF}\left( d_{\tilde{\lambda}}\tilde{\lambda}\right)=\psi$ and $\mathrm{LF}\left( d_{\tilde{\lambda}}\tilde{\lambda}\right)=\beta_0>0$, which implies that $d_{\tilde{\lambda}}\tilde{\lambda}\succ 1$ as desired.
\end{demonstration}

%----------------------------------------------------------------

\section{Pre-logarithms and derivations induced by decreasing automorphisms.}\label{section:decr-autom}

\subsection{Decreasing automorphisms and monomial series morphisms.}
\begin{defi}\label{defi:decreasing_endom}
 Let $(\Phi,\preccurlyeq)$ be a chain. A \textbf{decreasing endomorphism} $\sigma$ of $\Phi$
 is an order preserving map $\sigma : \Phi\rightarrow\Phi$, such that for all
 $\phi\in \Phi$, $\sigma(\phi)\prec\phi$. If this map is surjective, we call it a \textbf{decreasing automorphism}.
\end{defi}

\begin{remark}
Note that, if $\Phi$ has a decreasing endomorphism, then it has necessarily no least element. It would be interesting to characterize linear orderings which admit a decreasing endomorphism.
\end{remark}

%\begin{remark}
%In the definition of a decreasing endomorphism, $\Phi$ has necessarily no least element. In the case where $\Phi$ has a least element, say $\phi_m$, one can always define a decreasing map $\sigma:\Phi\setminus\{\phi_m\}\rightarrow\Phi$. Then one can \textbf{saturate} $\Phi$ thus :\\
%$\bullet$ consider some symbols $\sigma^{n}(\phi_m)$, $n\in\mathbb{N}^*$, ordered by: $\sigma^{k}(\phi_m)\prec\sigma^{l}(\phi_m)\Leftrightarrow k>l$;\\
%$\bullet$ consider $\tilde{\Phi}:=\Phi\cup \{\sigma^{n}(\phi_m),\ n\in\mathbb{N}^*\}$, as a new chain of fundamental monomials;\\
%$\bullet$ consider a group $\tilde{\Gamma}$ with $\mathrm{\textbf{H}}_{\mathrm{fin}}(\tilde{\Phi})\subset \tilde{\Gamma}\subset \mathrm{\textbf{H}}(\tilde{\Phi})$, and such that $\Gamma\subset \tilde{\Gamma}$.
%\sn Then the generalized series field $\tilde{\mathds{K}}:=\mathbb{R}((\tilde{\Gamma}))$ is an \emph{ordered field extension} of $\mathds{K}$.
%\sn
%Similarly, we leave it to the reader to verify that any decreasing endomorphism can be extended to a decreasing automorphism by adding new elements to the chain of fundamental monomials $\Phi$.
%\end{remark}

\begin{defi}\label{defi:monom-prelog}
 A pre-logarithm on $\mathds{K}$ is \textbf{monomial} if its restriction to the fundamental monomials has its image in the monomials: \begin{center}
$l:\Phi\rightarrow \mathbb{R}^*.\Gamma$.
\end{center}
\end{defi}

In \cite[Proposition 5.2]{matu-kuhlm:hardy-deriv-gener-series}, we study derivations on $\mathds{K}$ that are also called \textbf{monomial} (i.e. such that their restrictions to the fundamental monomials have their image in the monomials), and we prove that:
\begin{propo}\label{prop:monomial_case}
 A monomial derivation $d$ extends to a Hardy type series derivation on $\mathds{K}$ if and only if the condition (H3') holds.
\end{propo}
Here we prove that:

\begin{propo}\label{propo:monomial-prelog}
Let $d$ be a monomial Hardy type series derivation on $\mathds{K}$.
Assume that the set $\Theta=\left\{\theta^{(\phi)},\
\phi\in\Phi\right\}$ has no least element. Then there exists a
unique pre-logarithmic section $l$ on $\mathds{K}$ which is a series
morphism, for which the induced pre-logarithm is compatible with the
derivation. Moreover, this pre-logarithm verifies (GA).
\end{propo}
\begin{demonstration}
We just need to check the hypothesis of Theorem
\ref{theo:existence-prelog}. Indeed, for any $\phi$,
$\theta^{(\phi)}\neq \hat{\theta}$, which implies that assumption 1.
of Theorem \ref{theo:existence-prelog} holds. We compute now:
$$\mathrm{a.i.}(\theta^{(\phi)})= \displaystyle\frac{\theta^{(\phi)}}{\mathrm{LE}\left(\displaystyle\frac{ \theta^{(\phi)}}{\theta^{(\psi)}}\right)\mathrm{LT} \left(\displaystyle\frac{\psi '}{\psi}\right)}=\displaystyle\frac{1}{\mathrm{LE}\left(\displaystyle\frac{ \theta^{(\phi)}}{\theta^{(\psi)}}\right) \mathrm{LC}(\theta^{(\psi)})}
\displaystyle\frac{\theta^{(\phi)}}{\theta^{(\psi)}}$$ with $\mbox{ LF
}\left(\displaystyle\frac{\theta^{(\phi)}}{\theta^{(\psi)}}\right)=\psi$ (as in Lemma \ref{lemme:psi} with $\alpha=\theta^{(\phi)}$). Since $d$ is a Hardy
type derivation, by (H3') we have: $\mbox{ LF
}\left(\displaystyle\frac{\theta^{(\phi)}}{\theta^{(\psi)}}\right)
\prec\max\{\phi,\psi\}$ for any $\phi\neq\psi$. Consequently, $\phi=\max\{\phi,\psi\}\succ\psi$, which implies also that
$\theta^{(\phi)}\succ\theta^{(\psi)}$. Assumption 2 of Theorem
\ref{theo:existence-prelog} holds, as desired.
\end{demonstration}

\begin{ex}\label{ex:tressl}
We define the basic prelogarithmic section on $\mathds{K}$ by :
\[ l(\prod_{\phi \in \Phi} \phi^{\gamma_\phi}) = \sum_{\phi \in \Phi} \gamma_\phi \phi\]
Here (SL) is readily verified. The basic prelogarithmic section $l$
does {\it not} satisfy (GA) (e. g. $l(\phi)=\phi$ ). To remedy to
this problem, we fix a decreasing endomorphism \[\sigma\>: \Phi
\rightarrow \Phi.\] We define the \textbf{prelogarithmic section}
$l_{\sigma}$ \textbf{induced by }$\sigma$ as follows:
\[l_{\sigma}(\prod_{\phi \in \Phi} \phi^{\gamma_\phi}) = \sum_{\phi \in \Phi} \gamma_\phi \sigma(\phi)\>.\] The
\textbf{induced prelogarithm} (given in Definition \ref{defi:prelog}) is denoted by $l_{\sigma}$. We leave it to the reader to verify that $l_{\sigma}$ satisfies (GA) (see \cite[Ch. 3]{kuhl:ord-exp} for more details).
\mn
As an elementary but important illustration, take the following chain of infinitely increasing real germs at infinity (applying the usual comparison relations of germs) : \begin{center}
$\Phi:=\{\exp^n(x)\ ;\ n\in\mathbb{Z}\}$
\end{center}
where $\exp^n$ denotes for positive $n$, the $n$'th iteration of the real exponential function, for negative $n$, the $|n|$'s iteration of the logarithmic function, and for $n=0$ the identical map. The restriction of the (germ of the) natural logarithmic function $\log$ to $\Phi$ is such an embedding $\sigma$. We leave it to the reader to verify that its lifting as a pre-logarithm on $\mathds{K}$, extends the logarithmic function on the rational functions field $\mathbb{R}(\exp^n(x),\ n\in\mathbb{Z})$.
\end{ex}

\subsection{Defining a compatible monomial derivation from a series morphism.}

We study now the converse situation of Proposition
\ref{propo:monomial-prelog}. We consider the chain
$(\Phi,\preccurlyeq)$ endowed with a decreasing automorphism
$\sigma$, and the induced pre-logarithm $l_\sigma$. We want to know
when we can define a log-compatible Hardy type series derivation on
$\mathds{K}$, and describe it.

\begin{defi}\label{defi:orbit}
Given an ordered chain $(\Phi,\preccurlyeq)$, an element $\phi\in\Phi$ and an decreasing endomorphism $\sigma : \Phi\mapsto\Phi$, we call:\\
$\bullet$ the\textbf{ $\mathbb{Z}$-orbit} of $\phi$: $\mathcal{O}(\phi)=\{\sigma^k(\phi)\ |\ k\in\mathbb{Z}\}$;\\
$\bullet$ the \textbf{convex orbit} of $\phi$: $\mathcal{C}(\phi)=\{\psi\in\Phi\ |\ \exists k\in\mathbb{N},\ \sigma^k(\phi)\preccurlyeq\psi\preccurlyeq \sigma^{-k}(\phi)\}.$\\
$\bullet$  For any $\alpha=
\displaystyle\prod_{\phi\in \textrm{supp}\
\alpha}\phi^{\alpha_{\phi}}\in\Gamma$, any $\psi\in \Phi$ and any binary relation $\mathcal{R}\in\{\prec,\preccurlyeq,\succ,\succcurlyeq\}$, we denote   $S_\psi=\{\phi\in\textrm{supp}\
\alpha\ |\ \phi\mathcal{R}\psi\}$, and define the corresponding \textbf{truncation} of $\alpha$ as $\textrm{Tr}_{\mathcal{R}\psi}(\alpha):=\displaystyle\prod_{\phi\in
S_\psi}\phi^{\gamma_{\phi}}$.
\end{defi}
 \begin{nota}
 Given a family $\mathcal{F}\subset\Phi$ of representatives of the convex orbits of $\Phi$, given $\phi\in\mathcal{F}$, we denote $\mathcal{S}_{\mathcal{F},\phi}:=\{\psi\in\Phi\ |\ \phi\preccurlyeq\psi\prec\sigma^{-1}(\phi)\}$, and $\mathcal{S}_{\mathcal{F}}:=\displaystyle\bigcup_{\phi\in\mathcal{F}} \mathcal{S}_{\mathcal{F},\phi}$.
\end{nota}

\begin{propo}\label{propo:Hardy-deriv-monomiale}
Let $\sigma$ be a decreasing automorphism on $\Phi$, and $l_\sigma$
the induced pre-logarithm. There exists a log-compatible monomial
Hardy type series derivation on $\mathds{K}$ if and only if there
exists a map $\Phi\rightarrow\Gamma$, $\phi\mapsto\theta^{(\phi)}$,
such that:
\begin{description}
\item[(M)] for any $\phi\prec\psi\in\Phi$,  $\textrm{Tr}_{\succcurlyeq \mathcal{C}_\psi}\left(\displaystyle\frac{\theta^{(\psi)}}{ \theta^{(\phi)}}\right)=\textrm{Tr}_{\succcurlyeq \mathcal{C}_\psi}\left(\displaystyle\prod_{j=1}^{\infty} \displaystyle\frac{\sigma^j (\psi)}{\sigma^j(\phi)}\right)$,
with in particular,\\ for any $k\in\mathbb{N}$, $\theta^{(\sigma^{k}(\phi))}= \displaystyle\frac{\theta^{(\phi)}}{\prod_{j=1}^k\sigma^j(\phi)}$.
\end{description}
Moreover, given a family $\mathcal{F}$ of representatives of the various convex orbits of $\Phi$, such a derivation $d$ is unique up to the definition of the corresponding map $\mathcal{S}_{\mathcal{F}}\rightarrow \mathbb{R}^*\cdot\Gamma$, $\psi\mapsto t_\psi\theta^{(\psi)}$ (for arbitrary  $t_\psi\in\mathbb{R}^*$). In particular, when $\Phi$ admits only one convex orbit, say $\mathcal{C}_\phi$, then $d$ is unique up to the definition of $\theta^{(\phi)}\in\Gamma$, and $t_\psi\in\mathbb{R}^*$ for $\psi\in\Phi$ with $\phi\preccurlyeq\psi\prec\sigma^{-1}(\phi)$. More precisely,  we have $\theta^{(\psi)}= \theta^{(\phi)} \displaystyle\prod_{k=1}^{\infty}\displaystyle\frac{\sigma^k (\psi)}{\sigma^k(\phi)}$.
\end{propo}
\begin{demonstration} By Proposition \ref{prop:monomial_case}, the existence of a monomial Hardy type series derivation on $\mathds{K}$ reduces to the existence of a map $d:\Phi\rightarrow \mathbb{R}^*\Gamma$ such that (H3') holds.  By Proposition \ref{propo:compat}, such a series derivation is log-compatible if and only if (HL4) holds, which means, in the monomial case, that for any $\phi\in \Phi$, $(\sigma(\phi))'=\displaystyle\frac{\phi '}{\phi}= t_\phi.\theta^{(\phi)}$. But, $(\sigma(\phi))'=t_{\sigma(\phi)}.\theta^{(\sigma(\phi))}\sigma(\phi)$ by definition. Therefore, we obtain $t_{\sigma(\phi)}.\theta^{(\sigma(\phi))}=t_\phi. \displaystyle\frac{\theta^{(\phi)}}{\sigma(\phi)}$, and by induction, for any $k\in\mathbb{N}^*$ , $t_{\sigma^k(\phi)}.\theta^{(\sigma^k(\phi))}=t_\phi. \displaystyle\frac{\theta^{(\phi)}}{ \prod_{j=1}^k\sigma^j(\phi)}$, and $t_{\sigma^{-k}(\phi)}.\theta^{(\sigma^{-k}(\phi))}= t_\phi.\theta^{(\phi)} \displaystyle\prod_{j=0}^{k-1}\sigma^{-j}(\phi)$
 . Now, consider $\psi\in\Phi$ such that $\phi\preccurlyeq\psi\prec\sigma^{-1}(\phi)$, so $\sigma^k(\phi)\preccurlyeq \sigma^k(\psi)\prec \sigma^{k-1}(\phi)$ for any $k\in\mathbb{N}$. We deduce that $\displaystyle\frac{\theta^{(\phi)}}{\prod_{j=1}^k  \sigma^j(\phi)}\preccurlyeq \displaystyle\frac{\theta^{(\psi)}}{\prod_{j=1}^k \sigma^j(\psi)} \prec \displaystyle\frac{\theta^{(\phi)}}{\prod_{j=1}^{k-1} \sigma^j(\phi)}$, and equivalently $1\preccurlyeq\displaystyle\frac{\theta^{(\psi)}}{\theta^{(\phi)}} \displaystyle\prod_{j=1}^k\displaystyle\frac{\sigma^j(\phi)}{\sigma^j(\psi)} \prec \sigma^k(\phi)$. By letting $k$ tends to $+\infty$, we deduce that $1\preccurlyeq\displaystyle\frac{\theta^{(\psi)}}{\theta^{(\phi)}} \displaystyle\prod_{j=1}^{+\infty}\displaystyle\frac{\sigma^j(\phi)}{ \sigma^j(\psi)} \prec\chi$ for all $\chi\in\mathcal{C}_\phi$. For $\psi\in\Phi$ such that $\sigma^{-k}(\phi)\preccurlyeq\psi\prec\sigma^{-k-1}(\phi)$, we set $\tilde{\psi}:=\sigma^k(\psi)$. Then, 
 $\displaystyle\frac{\theta^{(\psi)}}{\theta^{(\phi)}} =\displaystyle\frac{\theta^{(\psi)}}{\theta^{(\tilde{\psi})}} \displaystyle\frac{\theta^{(\tilde{\psi})}}{\theta^{(\phi)}} =\displaystyle\prod_{j=1}^k\sigma^j(\psi) \displaystyle\frac{\theta^{(\tilde{\psi})}}{\theta^{(\phi)}}$. We are reduced to  the preceding case. Finally, assume that $\mathcal{C}_\phi\prec\mathcal{C}_\psi$,
 i.e. for any $k,l\in\mathbb{N}$, $\sigma^{-k}(\phi)\prec\sigma^l(\psi)$.
  By (H3'), we have $\mathrm{LF}\left( \displaystyle\frac{\theta^{(\sigma^l(\psi))}}{\theta^{(\phi)}}\right)=
   \mathrm{LF}\left( \displaystyle\frac{\theta^{(\psi)}}{\theta^{(\phi)} \prod_{j=1}^l\sigma^j(\psi)}\right)
   \prec \sigma^l(\psi)$,
   which implies that $\mathrm{LF}\left( \displaystyle\frac{\theta^{(\psi)}}{\theta^{(\phi)} \prod_{j=1}^\infty\sigma^j(\psi)}\right)\prec \mathcal{C}_\psi$. To conclude, it suffices to note that, in the present case, $\textrm{Tr}_{\succcurlyeq \mathcal{C}_\psi}\left(\displaystyle\prod_{j=1}^{\infty} \displaystyle\frac{\sigma^j (\psi)}{\sigma^j(\phi)}\right) =\displaystyle\prod_{j=1}^\infty\sigma^j(\psi)$.
\sn
Conversely, suppose now that there is a map $\Phi\rightarrow\Gamma$, $\phi\mapsto\theta^{(\phi)}$,
such that Condition (M) holds. We set $t_\phi:=1$ for any $\phi\in\Phi$.
 It remains to verify that (H3') and (HL4) hold for such a map $d:\Phi \rightarrow \Gamma$,
 $\phi\mapsto \phi '=\theta^{(\phi)}\phi$. Condition (HL4) holds since, for any $\phi\in\Phi$,
 $\sigma(\phi)'=\theta^{(\sigma(\phi))}\sigma(\phi)= \displaystyle\frac{\theta^{(\phi)}}{\sigma(\phi)}\sigma(\phi)
  =\displaystyle\frac{\phi '}{\phi}$. For (H3'), we consider $\phi\prec\psi\in\Phi$, and deduce from (M) that:
  $\mathrm{LF}\left(\displaystyle\frac{\theta^{(\psi)}}{ \theta^{(\phi)}}\right)=\sigma(\psi)\prec \psi$,
  and $\mathrm{LE}\left(\displaystyle\frac{\theta^{(\psi)}}{ \theta^{(\phi)}}\right)=1>0$.
\mn
Concerning the second part of the statement of the proposition, we observe from the preceding proof that,
 whenever we fix $\displaystyle\frac{\phi '}{\phi}:=t_\phi.\theta^{(\phi)}$, this determines the values
 of $\displaystyle\frac{\psi '}{\psi}$ for any $\psi\in \mathcal{O}(\phi)$. Then note that
  $\mathcal{S}_{\mathcal{F},\phi}$ is a family of representatives of the $\mathbb{Z}$-orbits included
  in $\mathcal{C}(\phi)$. Therefore, $(\mathcal{S}_{\mathcal{F},\phi})_{\phi\in\mathcal{F}}$ is a partition
   of $\Phi$, and $\mathcal{S}_{\mathcal{F}}$ is a family of representatives of the $\mathbb{Z}$-orbits of $\Phi$.
\end{demonstration}

\subsection{Examples.}
\noindent 1. Our purpose is to illustrate Proposition
\ref{propo:Hardy-deriv-monomiale}, in particular when the chain
$\Phi=\{\phi_i\ |\ i\in\mathbb{Z}\}$ is isomorphic to $\mathbb{Z}$.
Let $n\in\mathbb{N}^*$ be given. We consider the corresponding
automorphism  $\sigma$ of $\Phi$ defined by $\phi_i\mapsto
\phi_{i-n}$. For instance, we set  $\theta^{(\phi_{0})}:=1$. In
order to build a log-compatible monomial Hardy type series
derivation on $\mathds{K}$, we have to set
$\theta^{(\sigma^{-k}(\phi_{0}))}:=\theta^{(\phi_{kn})}=
\displaystyle\prod_{l=0}^{k-1}\phi_{ln}$, and
$\theta^{(\sigma^{k}(\phi_{0}))}=\theta^{(\phi_{-kn})}:=
\displaystyle\frac{1}{ \prod_{l=1}^{k}\phi_{-ln}}$, for any
$k\in\mathbb{N}$. Furthermore, for any $j\in\{1,\ldots,n-1\}$, we
have to set $\theta^{(\phi_j)}:=\displaystyle\prod_{l=1}^{+\infty}
\displaystyle\frac{\phi_{j-ln}}{ \phi_{-ln}}$. Then, for any
$k\in\mathbb{N}$,
$\theta^{(\sigma^{-k}(\phi_{j}))}=\theta^{(\phi_{j+kn})}:=
\displaystyle\prod_{l=1}^{+\infty}\displaystyle\frac{\phi_{j-ln}}{\phi_{-ln}}
\displaystyle\prod_{l=0}^{k-1}\phi_{j+ln}=
\displaystyle\frac{\prod_{l=-k+1}^{+\infty}\phi_{j-ln}}{
\prod_{l=1}^{+\infty}\phi_{-ln}}$, and
$\theta^{(\sigma^{k}(\phi_{j}))}=\theta^{(\phi_{j-kn})}:=
\displaystyle\prod_{l=k+1}^{+\infty} \displaystyle\frac{\phi_{j-ln}
}{\phi_{-ln}}\displaystyle\frac{1 }{\prod_{l=1}^{k}\phi_{-ln}} =
\displaystyle\frac{\prod_{l=k+1}^{+\infty}\phi_{j-ln}
}{\prod_{l=1}^{+\infty}\phi_{-ln}}$.\sn As an illustration with
germs of real functions at $+\infty$, consider for any
$i\in\mathbb{Z}$, $\phi_{2i}:=\log^{-i+1}(x)$ (with $\log^0(x):=x$),
and $\phi_{2i+1}:=\log^{-i+1}\circ g(x)$, where $g$ is a (ultimately
positive and differentiable) half compositional iterate of $\exp$
(i.e. $g\circ g(x)=\exp(x)$: see \cite[Section 6]{bosher:hardy}).
The automorphism of the chain $\Phi$ is the usual real logarithmic
function. We have: $\sigma(\phi_i)=\phi_{i+2}$.  By applying the
usual derivation with respect to $x$, for any $k\in\mathbb{N}^*$, we
compute: $\displaystyle\frac{\phi_{2k}'}{\phi_{2k}} =
\exp(x)\exp^2(x)\cdots\exp^{k-1}(x)=
\displaystyle\prod_{l=0}^{k-1}\phi_{2l}$, and
$\displaystyle\frac{\phi_{-2k}'}{\phi_{-2k}} =
\displaystyle\frac{1}{\log^k(x)\cdots\log(x)x}=
\displaystyle\prod_{l=1}^{k}\phi_{-l}$. Concerning the fundamental
monomial with odd indexes, following  Proposition
\ref{propo:Hardy-deriv-monomiale}, we have to set:
$\displaystyle\frac{\phi_{2k+1}'}{\phi_{2k+1}} :=
\displaystyle\frac{\prod_{l=-k+1}^{+\infty}\phi_{1-2l}}{
\prod_{l=1}^{+\infty}\phi_{-2l}}=
\displaystyle\frac{\prod_{l=-k+1}^{+\infty}\log^{l+1}\circ g(x)}{
\prod_{l=1}^{+\infty}\log^{l+1}}$, and
$\displaystyle\frac{\phi_{-2k+1}'}{\phi_{-2k+1}} :=
\displaystyle\frac{\prod_{l=k+1}^{+\infty}\phi_{1-2l}
}{\prod_{l=1}^{+\infty}\phi_{-2l}}=
\displaystyle\frac{\prod_{l=k+1}^{+\infty}\log^{l+1}\circ g(x)
}{\prod_{l=1}^{+\infty}\log^{l+1}}$. In particular,
$\displaystyle\frac{g'(x)}{g(x)}
=\displaystyle\frac{\phi_{-1}'}{\phi_{-1}} :=
\displaystyle\frac{\prod_{l=2}^{+\infty}\phi_{1-2l}
}{\prod_{l=1}^{+\infty}\phi_{-2l}}=
\displaystyle\prod_{l=1}^{+\infty}\displaystyle\frac{\log^{l+1}\circ
g(x) }{\log^{l}(x)}$. \sn It would be interesting to investigate the
possible analytic meaning of such a formal definition for the
derivative of $g$. \mn \noindent 2. The purpose now is to provide a
general example illustrating Proposition
\ref{propo:Hardy-deriv-monomiale}, with a uniform definition for the
$\theta^{(\phi)}$'s. Let $(\Phi,\preccurlyeq)$ be a chain endowed
with a decreasing automorphism $\sigma : \Phi\rightarrow\Phi$. Set
$\theta^{(\phi)}:=\displaystyle\prod_{k=1}^{+\infty}\sigma^k(\phi)$
and $t_\phi=1$ for any $\phi\in \Phi$. These momomials verify (M),
since for any $\phi\prec\psi\in\Phi$, we have
$\displaystyle\frac{\theta^{(\psi)}}{\theta^{(\phi)}}=
\displaystyle\prod_{k=1}^{+\infty}
\displaystyle\frac{\sigma^k(\psi)}{\sigma^k(\phi)}$. \sn In the case
of germs of real functions described at the end of Example \ref{ex:tressl} ($\Phi\approx\mathbb{Z}$), the present one can be seen as a
limit case. Indeed, instead of differentiating with respect to the
variable $x$, one may differentiate with respect to $\phi_i$, with
$i\rightarrow -\infty$. This can be viewed as a differentiation with
respect to a variable $\rho$ dominated by all the $\phi$'s in
$\Phi$: $\rho\prec\Phi$. In other words, differentiation with
respect to a \textit{translogarithm} (i.e. the compositional inverse
of a transexponential: see \cite{bosher:hardy}).

%-----------------------------------------------------------------------------

\section{Derivation on EL-series field}\label{sect:EL}
We consider $\mathds{K}$ endowed with a pre-logarithm $l$.
% which is a series morphism and verifies (GA), and with a log-compatible and strongly linear series derivation $d$.
The \textbf{exponential-logarithmic series} (EL-series for short)
field $(\mathds{K}^{\rm{ EL }},\log)$ corresponding to the
pre-logarithmic series field $(\mathds{K},l)$, is built as an
infinite towering extension of $\mathds{K}$, namely its exponential
closure (see below, \cite{kuhl:ord-exp} and \cite{kuhl-tressl:ELLE}
for details). Given a log-compatible series derivation $d$ on
$\mathds{K}$, the purpose of this section is to show how to extend
$d$ to a log-compatible series derivation (also denoted $d$) on
$\mathds{K}^{\rm{ EL }}$. If we assume moreover that $d$ is of Hardy
type, then so will be its extension.

\subsection{The exponential closure of a pre-logarithmic series field.}
Recall that the pre-logarithmic section $\ l\ :\ \Gamma\rightarrow \mathds{K}^{\succ 1}$ is an embedding of ordered groups. We denote by $\hat{\Gamma}=\mathds{K}^{\succ 1}\setminus l(\Gamma)$ the set complement of $l(\Gamma)$ in $\mathds{K}^{\succ 1}$, and by $\tilde{\Gamma}=e^{\hat{\Gamma}}$  a multiplicative copy of it (the choice of $e$ as abstract variable will result obvious from the definition of the new pre-logarithm $l^\sharp$ below). We endow the later with an ordering $\preccurlyeq$: $\forall\ e^a,e^b\in\tilde{\Gamma},\ e^a\prec e^b\Leftrightarrow a< b$. Then we define a new group $\Gamma^\sharp=\Gamma\cup\tilde{\Gamma}$ with the following
 multiplicative rule: if $\alpha^\sharp,\beta^\sharp\in\Gamma^\sharp$ both belong to $\Gamma$, multiply them there; similarly if they both belong to $\tilde{\Gamma}$. If $\alpha^\sharp=\alpha\in\Gamma$ and $\beta^\sharp=e^a\in\tilde{\Gamma}$ (i.e. $a\in\hat{\Gamma}$), then set $\alpha^\sharp.\beta^\sharp:=e^{l(\alpha)+a}$. Therefore $\Gamma^\sharp$ is a group extension of $\Gamma$.
\sn
We extend also $l$ to the following isomorphism:
\begin{center}
$\begin{array}{llcl}
l^\sharp\ :\ &(\Gamma^\sharp,.)&\rightarrow&(\mathds{K}^{\succ 1},+)\\
&\alpha^\sharp&\mapsto&l^\sharp(\alpha^\sharp)
\end{array}$
\end{center}
by defining $l^\sharp_{|\Gamma}:=l$, and for any $\alpha^\sharp=e^a\in\tilde{\Gamma}$, $l^\sharp(\alpha^\sharp):=l^\sharp(e^a)=a$.
Subsequently, we endow $\Gamma^\sharp$ with the ordering $\preccurlyeq$ defined as the transfer of the ordering $\leq$ on $\mathds{K}^{\succ 1}$. Hence it extends the ordering $\preccurlyeq$ on $\Gamma$.
\sn
We set $\mathds{K}^\sharp:=\mathbb{R}((\Gamma^\sharp))$, and the corresponding $(\mathds{K}^\sharp)^{\prec 1}$, $(\mathds{K}^\sharp)^{\preccurlyeq 1}$, $(\mathds{K}^\sharp)^{\succ 1}$ as before. Note that $\mathds{K}^{\succ 1}\subset(\mathds{K}^\sharp)^{\succ 1}$, so $l^\sharp\ :\ \Gamma^\sharp\rightarrow(\mathds{K}^\sharp)^{\succ 1}$ is a pre-logarithmic section. We extend it to a pre-logarithm $l^\sharp$ on $\mathds{K}^\sharp$ as in Definition \ref{defi:prelog}.
\mn
Repeating this process, we obtain inductively the $n^{th}$ extension of $(\mathds{K},l)$, denoted by $(\mathds{K}^{\sharp n},l^{\sharp n})$, $n\in\mathbb{N}$.
The corresponding EL-series field is defined as follows:
\begin{defi}
Set $\mathds{K}^{\rm{ EL
}}=\displaystyle\bigcup_{n\in\mathbb{N}}\mathds{K}^{\sharp n}$ and
$\log=\displaystyle\bigcup_{n\in\mathbb{N}}l^{\sharp n}$. We call
$(\mathds{K}^{\rm{ EL }},\log)$ the EL-series field over the
pre-logarithmic field $(\mathds{K},l)$.
\end{defi}
Note that $\log : ((\mathds{K}^{\rm{ EL }})^{>0}, \cdot)\rightarrow
(\mathds{K}^{\rm{ EL }}, +)$ is then an \textit{order preserving
isomorphism}. We denote by $\exp=\log^{-1}$ its inverse map.

\subsection{Extending derivations to the exponential closure.}
Consider a strongly linear (i.e. which verifies (D2)) and
log-compatible derivation $d$ on $\mathds{K}$. We show how to extend
$d$ to the corresponding EL-series field $\mathds{K}^{\rm{ EL }}$.
Note that this has been considered for fields of transseries in
\cite[Ch. 4.1.4]{Schm01}.
 However, our pre-logarithmic field $(\mathds{K},l)$ does not necessarily satisfy Axiom (T4)
  of \cite[Definition 2.2.1]{Schm01}.
\begin{theo}\label{theo:derivEL}
The strongly linear and log-compatible derivation $d$ on
$\mathds{K}$ extends to a strongly linear and log-compatible
derivation on $\mathds{K}^{\rm{ EL }}$, and this extension is
uniquely determined. Moreover, if $d$ is of Hardy type, then so is
its extension to $\mathds{K}^{\rm{ EL }}$.
\end{theo}

To prove the theorem, we proceed by induction along the towering extension process. \emph{Hence $(\mathds{K},l)$ represents from now until the end of this section, for simplicity of the notations, $(\mathds{K}^{\sharp n},l^{\sharp n})$ for some $n\in\mathbb{N}$}. We suppose $\mathds{K}$ endowed with a strongly linear and log-compatible derivation $d$, and require that its extension to $\mathds{K}^\sharp$ (also denoted by $d$) is also strongly linear and log-compatible:

\begin{lemma}\label{lemme:extension} For any $a^\sharp=\displaystyle\sum_{\alpha^\sharp\in\textrm{Supp}\ a^\sharp}a_{\alpha^\sharp}\alpha^\sharp\in\mathds{K}^\sharp$, if we set
\begin{center}
$d(a^\sharp)=(a^\sharp)'=\displaystyle\sum_{\alpha^\sharp\in\textrm{Supp}\ a^\sharp}a_{\alpha^\sharp}\alpha^\sharp (l^\sharp(\alpha^\sharp))'$,
\end{center}
 then $d$ is well-defined. Moreover, $d$ is the unique strongly linear and log-compatible derivation on $\mathds{K}^\sharp$ that extends $d$.\sn
Furthermore, if $d$ is a Hardy type derivation on $\mathds{K}$, then so it is on $\mathds{K}^\sharp$.
\end{lemma}
\begin{demonstration}
Consider $a^\sharp=\displaystyle\sum_{\alpha^\sharp\in\textrm{Supp}\ a^\sharp}a_{\alpha^\sharp}^\sharp\alpha^\sharp\in\mathds{K}^\sharp$. For any $\alpha^\sharp\in\textrm{Supp}\ a^\sharp$, we denote $\alpha^\sharp=\alpha$ if $\alpha^\sharp\in\Gamma$, and $\alpha^\sharp=e^a$ for some $a\in\hat{\Gamma}$ if $\alpha^\sharp\in\tilde{\Gamma}$. Then, by definition, we have:\begin{center}
$(a^\sharp)'=\displaystyle\sum_{\alpha\in(\textrm{Supp}\ a^\sharp)\cap\Gamma}a_{\alpha}^\sharp\alpha '+\displaystyle\sum_{e^a\in(\textrm{Supp}\ a^\sharp)\cap\tilde{\Gamma}}a_{e^a}^\sharp a'e^a$.
\end{center}
If we denote $a=\displaystyle\sum_{\alpha\in\textrm{Supp}\ a}a_\alpha\alpha\in\mathds{K}^{\succ 1}$, and $a'=\displaystyle\sum_{\beta\in\textrm{Supp}\ a'}b_{\beta}\beta\in\mathds{K}$, then:\begin{center}
$(a^\sharp)'=\displaystyle\sum_{\alpha\in(\textrm{Supp}\ a^\sharp)\cap\Gamma}a_{\alpha}^\sharp\alpha '+\displaystyle\sum_{e^a\in(\textrm{Supp}\ a^\sharp)\cap\tilde{\Gamma}}\displaystyle\sum_{\beta\in\textrm{Supp}\ a'} a_{e^a}^\sharp b_{\beta}e^{a+l(\beta)}$.
\end{center}
First, we verify that $(a^\sharp)'$ is well-defined. We set $S:=(\textrm{Supp}\ a^\sharp)\cap\Gamma$, and $\tilde{S}:=(\textrm{Supp}\ a^\sharp)\cap\tilde{\Gamma}$. Observe that $S$ and $\tilde{S}$ are anti-well-ordered subsets of $\Gamma$ and $\tilde{\Gamma}$ respectively. Hence, if we set $\hat{S}:=l^\sharp(\tilde{S})$, then $\hat{S}$ is anti-well-ordered in $(\mathds{K}^{\succ 1},\leq)$. The first sum is the derivative of $\displaystyle\sum_{\alpha\in S}a_{\alpha}^\sharp\alpha$, which is an element of $\mathds{K}$. By the induction hypothesis, it is well-defined. For the second sum, we have to show that the family $(a'e^a)_{a\in\hat{S}}$ is summable (see Definition \ref{defi:summ_fam}). As noted above, the elements of the support of this family are of the form $e^{a+l(\beta)}$, where $a\in\hat{S}$ and $\beta\in\textrm{Supp}\ a'$.
% its support $\displaystyle\bigcup_{e^a\in\tilde{S}} \displaystyle\bigcup_{\beta\in\textrm{Supp}\ a'}\{e^{a+l(\beta)}\}$ is anti-well-ordered in $\tilde{\Gamma}$, and that there can not be infinitely many terms contributing to a same monomial.
% By means of the automorphism $l^\sharp\ :\ \tilde{\Gamma}\rightarrow\hat{\Gamma}$, it is equivalent to show that the set $R:=\displaystyle\bigcup_{a\in\hat{S}} \displaystyle\bigcup_{\beta_a\in\textrm{Supp}\ a'}\{a+l(\beta_a)\}\subset\mathds{K}^{\succ 1}$ contains neither an infinite increasing sequence, nor a stationary one coming from differents $a+l(\beta_a)$'s.
Hence, to proceed by contradiction, we suppose that there is an increasing sequence $c_0\leq c_1\leq c_2\leq \cdots$ of elements of $\mathds{K}$ with $c_n:=a_n+l(\beta^{(n)})$, $a_n\in\hat{S}$, $\beta^{(n)}\in\textrm{Supp}\ a_n'$  for any $n$. Consider the corresponding sequence $(a_n)_{n\in\mathbb{N}}$. Since $\hat{S}$ is anti-well-ordered in $\mathds{K}^{\succ 1}$, it cannot have an increasing subsequence. Moreover, if it had a stationary subsequence, we would have a corresponding increasing subsequence of $l(\beta^{(n)})$'s. Since $l$ is order preserving, we would have an increasing subsequence of $\beta_n$'s, all of them belonging to the support of a same $a_{n_0}'$, contradicting the induction hypothesis. Therefore, by Lemma \ref{lemme:suite}, there is a strictly decreasing subsequence $a_{i_0}>a_{i_1}>a_{i_2}>\cdots$. Since the corresponding sequence  $c_{i_0}<c_{i_1}<\cdots$ is increasing, we must have $l(\beta^{(i_0)})<l(\beta^{(i_1)})<\cdots$, and equivalently $\beta^
 {(i_0)}\prec\beta^{(i_1)}\prec\cdots$. Subsequently, we observe that:
\begin{equation}\label{inequ1}\forall k<l\in\mathbb{N},\ \
0<a_{i_k}-a_{i_l}<l(\beta^{(i_l)})-l(\beta^{(i_k)}) =l\left(\displaystyle\frac{\beta^{(i_l)}}{\beta^{(i_k)}}\right)
\end{equation}
with $\beta^{(i_k)}\in\textrm{Supp}\ a_{i_k}'$ and $\beta^{(i_l)}\in\textrm{Supp}\ a_{i_l}'$. By the induction hypothesis, we denote $\beta^{(i_n)}=\alpha^{(i_n)}\gamma^{(i_n)}$, where $\alpha^{(i_n)}\in\textrm{Supp}\ a_{i_n}$ and $\gamma^{(i_n)}\in\textrm{Supp}\ \displaystyle\frac{(\alpha^{(i_n)})'}{\alpha^{(i_n)}}=\textrm{Supp}\ l(\alpha^{(i_n)})'$ for any $n\in\mathbb{N}$. We observe that $\alpha^{(i_n)}\in\Gamma^{\succ 1}$ since $a_{i_n}\in\mathds{K}^{\succ 1}$, and  $l\left(\displaystyle\frac{\beta^{(i_l)}}{\beta^{(i_k)}}\right)= l\left(\displaystyle\frac{\alpha^{(i_l)}}{\alpha^{(i_k)}}\right)+ l\left(\displaystyle\frac{\gamma^{(i_l)}}{\gamma^{(i_k)}}\right)$ for any $k,l$. Consider the sequence $(\alpha^{(i_n)})_{n\in\mathbb{N}}$. If it had a strictly decreasing subsequence, we could define the series $\tilde{a}=\displaystyle\sum_{n\in\mathbb{N}}\alpha^{(i_n)}$. But the corresponding increasing subsequence of $(\beta^{(i_n)})_{n\in\mathbb{N}}$ would be included in $\textrm{Supp}\ \tilde{a}'$, contradiction. Neither can $(\alpha^{(i_n)})_{n\in\mathbb{N}}$ have any stationary subsequence. Indeed, the corresponding increasing subsequence of $(\beta^{(i_n)})_{n\in\mathbb{N}}$ would be included in $\textrm{Supp}\ (\alpha^{(i_{n_0})})'$ for a fixed $n_0$. Hence, by Lemma \ref{lemme:suite}, $(\alpha^{(i_n)})_{n\in\mathbb{N}}$ has a strictly increasing subsequence, say $(\alpha^{(j_n)})_{n\in\mathbb{N}}$.
%  and the corresponding increasing ones $(\beta^{(j_n)})_{n\in\mathbb{N}}$, $(c_{j_n})_{n\in\mathbb{N}}$, and decreasing one $(a_{j_n})_{n\in\mathbb{N}}$.
Observe now that:
\begin{equation}\label{inequ2}\forall k<l\in\mathbb{N},\ \
1\prec\displaystyle\frac{\gamma^{(j_k)}}{\gamma^{(j_l)}}\preccurlyeq \displaystyle\frac{\alpha^{(j_l)}}{\alpha^{(j_k)}}
\end{equation}
It implies that $0<l\left(\displaystyle\frac{\gamma^{(j_k)}}{\gamma^{(j_l)}}\right) \leq l\left(\displaystyle\frac{\alpha^{(j_l)}}{\alpha^{(j_k)}}\right)<l(\alpha^{(j_l)})$ (since $\alpha^{(j_k)}\in\Gamma^{\succ 1}$). But, since $l$ verifies (GA), we also have $l(\alpha^{(j_l)})\prec\alpha^{(j_l)}$. Therefore, $l\left(\displaystyle\frac{\alpha^{(j_l)}}{\alpha^{(j_k)}}\right) \prec \alpha^{(j_l)}$, which implies that $l\left(\displaystyle\frac{\gamma^{(j_k)}}{\gamma^{(j_l)}}\right) \prec\alpha^{(j_l)}$. We obtain that $l\left(\displaystyle\frac{\beta^{(i_l)}}{\beta^{(i_k)}}\right) \prec\alpha^{(j_l)}$. But, we deduce from (\ref{inequ1}) that $a_{i_k}-a_{i_l}\prec\alpha^{(j_l)}$. It implies that the term with monomial $\alpha^{(j_l)}$ in $a_{i_l}$ has been cancelled by a term in $a_{i_k}$. Therefore, $\alpha^{(j_l)}\in\textrm{Supp}\ a_{i_k}$, so $\beta^{(j_l)}\in\textrm{Supp}\ a_{i_k}'$. By a straightforward induction, we obtain that the strictly increasing sequence $(\beta^{(j
 _n)})_{n\in\mathbb{N}}$ is included in $\textrm{Supp}\ a_{i_0}'$, contradiction. The extension of $d$ to $\mathds{K}^\sharp$ is well-defined.
\sn
The proofs that such an extension is a log-compatible derivation follow straightly from the definitions and are left to the reader.
\sn
Moreover, we observe that the extension of $d$ to $\mathds{K}^\sharp$ is uniquely determined by its definition, since we suppose that it is strongly linear and log-compatible: \begin{center}
$d(a^\sharp)=(a^\sharp)'=\displaystyle\sum_{\alpha^\sharp\in\textrm{Supp}\ a^\sharp}a_{\alpha^\sharp}^\sharp\alpha^\sharp (l^\sharp(\alpha^\sharp))'$.
\end{center}
Suppose now additionally that $d$ is a Hardy type derivation on $\mathds{K}$. To prove that $d$ verifies
 l'Hospital's rule (HD2) on $\mathds{K}^\sharp$, it suffices to prove it for its monomials. Hence, we consider $\alpha^\sharp,\beta^\sharp\in\Gamma^\sharp\backslash\{1\}$ with $\alpha^\sharp\prec\beta^\sharp$. It means that $l^\sharp(\alpha^\sharp)<l^\sharp(\beta^\sharp)$ in $\mathds{K}^{\succ 1}$, which is equivalent, by the induction hypothesis, to: $l^\sharp(\alpha^\sharp)'<l^\sharp(\beta^\sharp)'$. Therefore: $(\alpha^\sharp)'=\alpha^\sharp l^\sharp(\alpha^\sharp)'\prec \beta^\sharp l^\sharp(\beta^\sharp)'=(\beta^\sharp)'$.
\sn
To determine the subfield of constants of $\mathds{K}^\sharp$, suppose now that there exists $a^\sharp=\displaystyle\sum_{\alpha^\sharp\in\textrm{Supp}\ a^\sharp}a^\sharp_{\alpha^\sharp}\alpha^\sharp \in\mathds{K}^\sharp\setminus\mathbb{R}$ such that $(a^\sharp)'=0$. We denote as before:\begin{center}
$\begin{array}{lcl}
(a^\sharp)'&=&\displaystyle\sum_{\alpha^\sharp\in\textrm{Supp}\ a^\sharp}a_{\alpha^\sharp}^\sharp\alpha^\sharp (l^\sharp(\alpha^\sharp))'\\
&=&\displaystyle\sum_{\alpha\in S}a^\sharp_{\alpha}\alpha '+\displaystyle\sum_{e^a\in \tilde{S}}a^\sharp_{e^a}a'e^a\\
&=&\displaystyle\sum_{\alpha\in S}a^\sharp_{\alpha}\alpha '+\displaystyle\sum_{a\in\hat{S}}\displaystyle\sum_{\beta\in\textrm{Supp}\ a'} a^\sharp_{e^a}b_{\beta}e^{a+l(\beta)},
\end{array}$
\end{center}
where $S:=(\textrm{Supp}\ a^\sharp)\cap\Gamma$, $\tilde{S}:=(\textrm{Supp}\ a^\sharp)\cap\tilde{\Gamma}$, and $\hat{S}:=l^\sharp(\tilde{S})$.
Set $\alpha^\sharp_0:=\max(\textrm{Supp}\ a\backslash\{1\})$. There are two possibilities. Either $\alpha^\sharp_0=\alpha_0\in\Gamma$. By the induction hypothesis, it implies that there is $\beta_0\in\textrm{Supp}\ \alpha_0'$. So the corresponding term in $(a^\sharp)'$ must have been cancelled by the leading term of the second sum. Or $\alpha^\sharp_0=e^{a_0}$ for some $a_0\in\hat{\Gamma}$.
 Then there exists $\hat{\alpha}_0:=\textrm{LM}(a_0)\succ 1$ and
 $\hat{\beta}_0:=\textrm{LM}(\hat{\alpha}_0')\neq 0$ (by the induction hypothesis).
  So $e^{a_0+l(\hat{\beta}_0)}$ is the leading monomial of the second sum.
  The corresponding term in $(a^\sharp)'$ must have been cancelled by the leading term of the first sum,
  say $\beta_0$. Therefore, in the two cases, there must be an equality $e^{a_0+l(\hat{\beta}_0)}=\beta_0$,
   which is equivalent to: $\displaystyle\frac{e^{a_0+l(\hat{\beta}_0)}}{\beta_0}=1$.
   But  $\displaystyle\frac{e^{a_0+l(\hat{\beta}_0)}}{\beta_0}= e^{a_0+l(\hat{\beta}_0)-l(\beta_0)}$.
   So we should have $a_0+l(\hat{\beta}_0)-l(\beta_0)=a_0+ l\left(\displaystyle\frac{\hat{\beta}_0}{\beta_0}\right)
   =0$, which is absurd since $a_0\in\hat{\Gamma}=\mathds{K}^{\succ 1}\setminus l(\Gamma)$.
\sn
It remains to prove (HD3) on $\mathds{K}^\sharp$. We consider $a^\sharp,b^\sharp\in\mathds{K}_\sharp$ with $|a^\sharp|\succ |b^\sharp|\succ 1$. Note that, by replacing $a^\sharp$ and $b^\sharp$ by $-a^\sharp$, $\pm\displaystyle\frac{1}{a^\sharp}$, and $-b^\sharp$, $\pm\displaystyle\frac{1}{b^\sharp}$ respectively, we still have $|a^\sharp|\succ |b^\sharp|\succ 1$, and the leading monomials of $\displaystyle\frac{(a^\sharp) '}{a^\sharp}$ and $\displaystyle\frac{(b^\sharp) '}{b^\sharp}$ are preserved. So we can suppose without loss of generality that $a^\sharp\succ b^\sharp\succ 1$. Consequently, $l^\sharp(a^\sharp)>l^\sharp(b^\sharp)>0$ in $\mathds{K}$. This implies that $\displaystyle\frac{(a^\sharp) '}{a^\sharp} =l^\sharp(a^\sharp)'\succcurlyeq l^\sharp(b^\sharp)'=\displaystyle\frac{(b^\sharp) '}{b^\sharp}$. Moreover, $a^\sharp\succ\succ b^\sharp\succ 1$ if and only if $l^\sharp(a^\sharp)\succ l^\sharp(b^\sharp)$ in $\mathds{K}$, which means that $\displaystyle\frac{(a^\sharp) '}{a^\sharp}=l^\sharp(a^\sharp)'\succ l^\sharp(b^\sharp)'=\displaystyle\frac{(b^\sharp) '}{b^\sharp}$.\sn
This concludes the proof of Lemma \ref{lemme:extension}, and so the one of Theorem \ref{theo:derivEL}.
\end{demonstration}

%---------------------------------------------------------------------------

\section{Asymptotic integration and integration on EL-series}\label{sect:integr}
Let $(\mathds{K},l,d)$ be a pre-logarithmic series field endowed
with a strongly linear and log-compatible Hardy type
 derivation $d$. Recall that  $\hat{\theta}:= \mbox{ g.l.b. }_{\preccurlyeq}\left\{\theta^{(\phi)},\ \phi\in\Phi\right\}$,
  whenever it exists in $\Gamma$. In particular, in the case where $d$ is a series derivation,
   $\hat{\theta}= \mbox{ g.l.b. }_\preccurlyeq\left\{\textrm{LM}\left(\displaystyle\frac{a'}{a}\right);\ a\in\mathds{K}^*,\
   a\ \nasymp\ 1\right\}$
   (see Proposition \ref{propo:theta-hat}).

\begin{theo}\label{theo:integ-asymp-EL}
Let $(\mathds{K}^{\rm{ EL }},\log,d)$ be the induced differential
EL-series field as in Theorem \ref{theo:derivEL}.
 A series $a\in\mathds{K}^{\rm{ EL }}$ admits an asymptotic integral if and only if $a\nasymp \hat{\theta}$.
\end{theo}
\begin{demonstration} Recall that the induced derivation $d$ on $\mathds{K}^{\rm{ EL }}$ is itself strongly linear, log-compatible and of Hardy type (Theorem
\ref{theo:derivEL}). We proceed by induction along the towering
extension construction of $\mathds{K}^{\rm{ EL }}$. The initial step
is given by Theorem \ref{theo:integ-asymp-prelog}. Consider
$(\mathds{K},l,d)$ as in the statement of the theorem. We want to
show that its  extension $(\mathds{K}^\sharp,l^\sharp)$ verifies
 $\mbox{ g.l.b. }_\preccurlyeq\left\{\textrm{LM}\left(\displaystyle\frac{(a^\sharp) '}{a^\sharp}\right)
 ;\ a^\sharp\in\mathds{K}^\sharp\backslash\{0\},\ a^\sharp\ \nasymp\ 1\right\}=\hat{\theta}$
 (indeed, we will then obtain the desired result by applying Theorem \ref{theo:integ-asymp-prelog}
  to $\mathds{K}^\sharp=\mathbb{R}((\Gamma^\sharp))$). Let $1\nasymp a^\sharp\in\mathds{K}^\sharp\backslash\{0\}$.
   We denote $\alpha^\sharp:=\textrm{LM}(a^\sharp)$.
   By (HD2), $\textrm{LM}\left(\displaystyle\frac{(a^\sharp) '}{a^\sharp}\right)
   =\textrm{LM}\left(\displaystyle\frac{(\alpha^\sharp) '}{\alpha^\sharp}\right)$.
   There are two possibilities. Either $\alpha^\sharp=\alpha\in\Gamma$,
    which implies that $\textrm{LM}\left(\displaystyle\frac{(\alpha^\sharp) '}
    {\alpha^\sharp}\right)\succcurlyeq\hat{\theta}$. Or $\alpha^\sharp=e^a\in\tilde{\Gamma}$
    for some $a\in\hat{\Gamma}=\mathds{K}^{\succ 1}\backslash l(\Gamma)$.
    We denote $\alpha:=\textrm{LM}(a)$.
    Then $\textrm{LM}\left(\displaystyle\frac{(\alpha^\sharp) '}{\alpha^\sharp}\right)=
    \textrm{LM}(a')=\textrm{LM}(\alpha ')=\alpha \textrm{LM}\left(\displaystyle\frac{(\alpha) '}{\alpha}\right)$.
    Since $\textrm{LM}\left(\displaystyle\frac{(\alpha) '}{\alpha}\right)\succcurlyeq\hat{\theta}$,
     and $\alpha\succ 1$,
     then $\textrm{LM}\left(\displaystyle\frac{(\alpha^\sharp) '}{\alpha^\sharp}\right)\succ\theta$.
      Therefore, by induction we obtain that
      $\hat{\theta}=\mbox{ g.l.b. }_\preccurlyeq\left\{\textrm{LM}\left(\displaystyle\frac{(\hat{a}) '}{\hat{a}}\right)
      ;\ \hat{a}\in\mathds{K}^{\sharp n}\backslash\{0\},\ \hat{a}\ \nasymp\ 1\right\}$
      for any $n\in\mathbb{N}$, so $\hat{\theta}=\mbox{ g.l.b. }_\preccurlyeq\left
      \{\textrm{LM}\left(\displaystyle\frac{(\hat{a}) '}{\hat{a}}\right);\ \hat{a}\in\mathds{K}^{\rm{ EL }}\backslash\{0\},
      \ \hat{a}\ \nasymp\ 1\right\}$ as desired.
\end{demonstration}
Denote $\Gamma^{\rm{ EL }}:=\mathrm{LM}(\mathds{K}^{\rm{ EL
}})=\bigcup_{n\in\mathbb{N}} \Gamma^{\sharp n}$.

\begin{coro}
The EL-series field $\mathds{K}^{\rm{ EL }}$ is closed under
integration if and only if $\hat{\theta}\notin \Gamma^{\rm{ EL
}}$.
\end{coro}
\begin{demonstration}
By \cite[Theorem 55 and 56]{kuhl:Hensel-lemma}, it suffices to prove that
$\hat{\theta}\notin \Gamma^{\rm{ EL }}$ implies that any element of
$\mathds{K}^{\rm{ EL }}$ has an asymptotic integral, which was
proved in the preceding theorem.
\end{demonstration}

\bibliographystyle{amsalpha}
\def\cprime{$'$}
\providecommand{\bysame}{\leavevmode\hbox to3em{\hrulefill}\thinspace}
\providecommand{\MR}{\relax\ifhmode\unskip\space\fi MR }
\providecommand{\MRhref}[2]{%
  \href{http://www.ams.org/mathscinet-getitem?mr=#1}{#2}
}
\providecommand{\href}[2]{#2}

\end{document}